\documentclass[a4paper,12 pt]{article}
\usepackage{amsfonts}
\usepackage{amssymb}
\usepackage{latexsym}
\usepackage{amsmath}
\usepackage{amscd}
\usepackage{amsthm}
\usepackage{graphicx}
\usepackage{epsfig}
\usepackage{hyperref}
\usepackage{pst-all} 
\usepackage{psfrag}
\usepackage{url}

\newtheorem{theorem}{Theorem}[section]
\newtheorem{definition}[theorem]{Definition}

\newtheorem{proposition}[theorem]{Proposition}
\newtheorem{lemma}[theorem]{Lemma}
\newtheorem{remark}[theorem]{Remark}
\newtheorem{corollary}[theorem]{Corollary}

\newtheorem{prop}[theorem]{Proposition}

\newcommand{\E}{{\mathbb E}}
\newcommand{\N}{{\mathbb N}}
\newcommand{\Z}{{\mathbb Z}}
\renewcommand{\P}{{\mathbb P}}
\newcommand{\Q}{{\mathbb Q}}
\newcommand{\R}{{\mathbb R}}

\renewcommand{\u}{{\mathbf u}}

\newcommand{\z}{{\mathbf z}}
\newcommand{\TT}{{\mathcal T}}
\newcommand{\WW}{{\mathcal W}}

\usepackage[normalem]{ulem}
\usepackage{color}

\newcommand{\cg}{\color{green}}

\newcommand{\footremember}[2]{%
\footnote{#2}
\newcounter{#1}
\setcounter{#1}{\value{footnote}}%
}

\title{Continuum random tree as the scaling limit for a drainage network 
model: a Brownian web approach}
\author{%
  Kumarjit Saha\footremember{TIFR}{Ashoka University}
  \footnote{This work is benefited from the support of the AIRBUS 
  Group Corporate Foundation Chair in Mathematics of Complex Systems established in ICTS-TIFR
   and supported in part by ISF-UGC grant.}
  }

\begin{document}

\maketitle

\begin{abstract}
We consider  the tributary
structure of Howard’s drainage network model studied by Gangopadhyay et al. \cite{GRS04}. 
Conditional on the event that the tributary survives up to time $n$, we show that, as a sequence of random metric spaces, scaled tributary converges in distribution to a continuum random tree with respect to Gromov Hausdorff topology.
This verifies a prediction made by Aldous \cite{A93} for a simpler model (where  paths are independent till they coalesce) but for a different conditional set up. The limiting continuum random tree is slightly different from what was surmised in
\cite{A93}. Our proof uses the fact that there exists a dual process such that the original network and it's dual jointly converge to the Brownian web and its dual. The  limiting continuum random tree is universal in the sense that for all discrete drainage network models with non-crossing paths 
in the basin of attraction of the Brownian web, the conditional scaled tributaries converge to the  same scaling limit.  
\end{abstract}



\section{Introduction}

Trees are often used to represent genealogical structures,
communication network modelling, optimisation and
river basin study. During the last three decades there has been an interest in 
understanding the  scaling limits of various discrete random trees. 
Most notably in \cite{A91}, Aldous showed that an appropriately scaled Galton-Watson tree with a finite variance critical offspring distribution and conditioned to have  total population size $n$,  as $n \to \infty$, converges to a continuum random tree (CRT) famously known as the Brownian CRT. In the recent years, these results have been extended for more general continuum random trees (e.g., see \cite{LD}). 
These concepts and results have found 
many other applications also, e.g., understanding the scaling limits for thhe sequence of components 
of critical Erdos-Renyi random graph (\cite{A97}, \cite{ABG12}). Further this 
limiting behaviour appear to be universal in the sense that Aldous' limiting picture for  
critical random graph has since been extended to ``immigration" models of random graphs
\cite{AP00}, hypergraphs \cite{G05} and to random graphs with fixed degree \cite{NP10}.

In this paper we show that a sequence of appropriately scaled  discrete trees obtained from a drainage network model converge in distribution to a continuum random tree which is different from the Brownian CRT. 
Various statistical models of drainage networks have been proposed (see \cite{RR97} for a detailed survey). In this paper, we
study the `tributary' tree of a two-dimensional drainage network called
the Howard’s model of headward growth and branching \cite{H71}. 
In order to present our result, we first describe Howard's  network model. Fix $p \in (0,1)$ and let $\{B_{\u} : \u = (\u(1), \u(2)) \in \Z^2\}$ be an i.i.d. collection of Bernoulli random variables with success probability $p$. In what follows, a vertex $(x,t)$ with $B_{(x,t)} = 1$ will be called as an \textit{open} vertex.  Let $\{U_{\u}: \u \in \Z^2\}$ be another collection of random variables, independent of the collection $\{B_{\u} : \u \in \Z^2\}$    taking values $1$ and $-1$ with equal probabilities. For $(x,t) \in \Z^2$, let $k_0 := \min\{|k| : k \in \Z, B_{(x+k,t+1)} = 1\}$ and we define
\begin{align*}
h(x,t) := 
\begin{cases}
(x+k_0, t+1 ) & \text{ if }B_{(x - k_0, t+1 )} \notin V \\
(x- k_0, t+1 ) & \text{ if }B_{(x + k_0, t+1 )} \notin V \\
(x + U_{(x,t)} k_0, t+1 ) & \text{ otherwise.}
\end{cases}
\end{align*}
This gives a random graph ${\cal G}$ with vertex set $V := \{\u : B_\u = 1\}$, i.e., the set of open vertices, and
 edge set $E := \{\langle (x,t), h (x,t)\rangle : (x,t) \in V\}$. In other words, each open vertex connects to the nearest open vertex at the next level and in case of ties, chooses one of them uniformly independent of everything else. This gives a stochastic drainage network model, where each open vertex acts as a source vertex and from each open vertex $(x,t)$, water comes out and flows to $h(x,t)$ along the channel given by the edge $\langle (x,t), h(x,t)\rangle$. Clearly by construction, ${\cal G}$ does not have a cycle or loop almost surely.

Gangopadhyay et al. \cite{GRS04} studied this random graph in detail and showed the following: 
\begin{theorem}
\label{thm:GRSTreeForest}
${\mathcal G}$ is connected and there is no bi-infinite path in ${\mathcal G}$ almost surely.
\end{theorem}
It should be mentioned here that, Gangopadhyay et al. \cite{GRS04} studied this model for $\Z^d$ with $d\geq 2$ and observed 
a tree-forest dichotomy behaviour depending on dimensions. In this paper, we are interested about $d=2$ only and the construction that we mentioned here is different from that of \cite{GRS04} but 
generates the same process for $d = 2$. 

For $k \geq 1$ and for $(x,t)\in \Z^2$, let $h^k(x,t) = h(h^{k-1}(x,t))$ with $h^0(x,t) = (x,t)$.
For $(x, t) \in \Z^2$,
we now define the  cluster or the tributary at $(x,t)$ (see Figure \ref{fig:GRSDual}),
consisting of all the source vertices whose water flows through $(x,t)$, as
\begin{align}
\label{eq:ClusterLevel}
C(x,t) & := \{(y,s) \in \Z^2 : h^l(y,s) = (x,t) \text{ for some }l \geq 0\}.
\end{align}
The set of edges in cluster $C(x,t)$ is denoted by 
$$
E(x,t) := \{\langle (y,s), h(y,s) \rangle : (y,s) \in C(x,t), (y,s)\neq (x,t)\}.
$$

\begin{figure}[htb]
\begin{center}
\begin{pspicture}(0,0)(12,6)
\psgrid[subgriddiv=0,griddots=10,gridlabels=0]
\pscircle[fillcolor=black,fillstyle=solid](1,0){.1}
\pscircle[fillcolor=black,fillstyle=solid](1,1){.1}
\pscircle[fillcolor=black,fillstyle=solid](1,2){.1}
\pscircle[fillcolor=black,fillstyle=solid](1,3){.1}
\pscircle[fillcolor=black,fillstyle=solid](2,4){.1}
\pscircle[fillcolor=black,fillstyle=solid](0,5){.1}
\pscircle[fillcolor=black,fillstyle=solid](0,6){.1}

\psline[linewidth=1.5pt]{->}(1,0)(1,1)
\psline[linewidth=1.5pt]{->}(1,1)(1,2)
\psline[linewidth=1.5pt]{->}(1,2)(1,3)
\psline[linewidth=1.5pt]{->}(1,3)(2,4)
\psline[linewidth=1.5pt]{->}(2,4)(0,5)
\psline[linewidth=1.5pt]{->}(0,5)(0,6)

\pscircle[fillcolor=black,fillstyle=solid](3,1){.1}
\pscircle[fillcolor=black,fillstyle=solid](4,2){.1}
\pscircle[fillcolor=black,fillstyle=solid](4,3){.1}
\psline[linewidth=1.5pt]{->}(3,1)(4,2)
\psline[linewidth=1.5pt]{->}(4,2)(4,3)
\psline[linewidth=1.5pt]{->}(4,3)(5,4)

\pscircle[fillcolor=black,fillstyle=solid](4,5){.1}
\pscircle[fillcolor=black,fillstyle=solid](4,6){.1}
\psline[linewidth=1.5pt]{->}(4,5)(4,6)

\pscircle[fillcolor=black,fillstyle=solid](4,0){.1}
\pscircle[fillcolor=black,fillstyle=solid](5,1){.1}
\pscircle[fillcolor=black,fillstyle=solid](6,2){.1}
\pscircle[fillcolor=black,fillstyle=solid](6,3){.1}
\pscircle[fillcolor=black,fillstyle=solid](5,4){.1}
\pscircle[fillcolor=black,fillstyle=solid](5,5){.1}
\pscircle[fillcolor=black,fillstyle=solid](6,6){.1}

\psline[linewidth=1.5pt]{->}(4,0)(5,1)
\psline[linewidth=1.5pt]{->}(5,1)(6,2)
\psline[linewidth=1.5pt]{->}(6,2)(6,3)
\psline[linewidth=1.5pt]{->}(6,3)(5,4)
\psline[linewidth=1.5pt]{->}(5,4)(5,5)
\psline[linewidth=1.5pt]{->}(5,5)(6,6)

\pscircle[fillcolor=black,fillstyle=solid](9,1){.1}
\pscircle[fillcolor=black,fillstyle=solid](9,2){.1}
\pscircle[fillcolor=black,fillstyle=solid](8,3){.1}
\pscircle[fillcolor=black,fillstyle=solid](9,4){.1}
\pscircle[fillcolor=black,fillstyle=solid](8,5){.1}
\pscircle[fillcolor=black,fillstyle=solid](9,6){.1}

\pscircle[fillcolor=black,fillstyle=solid](7,4){.1}
\pscircle[fillcolor=black,fillstyle=solid](8,1){.1}
\pscircle[fillcolor=black,fillstyle=solid](8,2){.1}
\psline[linewidth=1.5pt]{->}(8,1)(8,2)
\psline[linewidth=1.5pt]{->}(8,2)(8,3)
\psline[linewidth=1.5pt]{->}(7,4)(8,5)
\psline[linewidth=1.5pt]{->}(9,1)(9,2)
\psline[linewidth=1.5pt]{->}(9,2)(8,3)
\psline[linewidth=1.5pt]{->}(8,3)(9,4)
\psline[linewidth=1.5pt]{->}(9,4)(8,5)
\psline[linewidth=1.5pt]{->}(8,5)(9,6)

\pscircle[fillcolor=black,fillstyle=solid](11,0){.1}
\pscircle[fillcolor=black,fillstyle=solid](12,1){.1}
\pscircle[fillcolor=black,fillstyle=solid](11,2){.1}
\pscircle[fillcolor=black,fillstyle=solid](12,3){.1}
\pscircle[fillcolor=black,fillstyle=solid](12,4){.1}
\pscircle[fillcolor=black,fillstyle=solid](11,5){.1}
\pscircle[fillcolor=black,fillstyle=solid](12,6){.1}
\pscircle[fillcolor=black,fillstyle=solid](10,4){.1}
\psline[linewidth=1.5pt]{->}(10,4)(11,5)
\psline[linewidth=1.5pt]{->}(11,0)(12,1)
\psline[linewidth=1.5pt]{->}(12,1)(11,2)
\psline[linewidth=1.5pt]{->}(11,2)(12,3)
\psline[linewidth=1.5pt]{->}(12,3)(12,4)
\psline[linewidth=1.5pt]{->}(12,4)(11,5)
\psline[linewidth=1.5pt]{->}(11,5)(12,6)

\pscircle[fillcolor=gray,fillstyle=solid](2.5,0){.1}
\pscircle[fillcolor=gray,fillstyle=solid](7.5,0){.1}
\pscircle[fillcolor=gray,fillstyle=solid](2,1){.1}
\pscircle[fillcolor=gray,fillstyle=solid](4,1){.1}
\pscircle[fillcolor=gray,fillstyle=solid](6.5,1){.1}
\pscircle[fillcolor=gray,fillstyle=solid](8.5,1){.1}
\pscircle[fillcolor=gray,fillstyle=solid](10.5,1){.1}

\pscircle[fillcolor=gray,fillstyle=solid](2.5,2){.1}
\pscircle[fillcolor=gray,fillstyle=solid](5,2){.1}
\pscircle[fillcolor=gray,fillstyle=solid](7,2){.1}
\pscircle[fillcolor=gray,fillstyle=solid](8.5,2){.1}
\pscircle[fillcolor=gray,fillstyle=solid](10,2){.1}

\pscircle[fillcolor=gray,fillstyle=solid](2.5,3){.1}
\pscircle[fillcolor=gray,fillstyle=solid](5,3){.1}
\pscircle[fillcolor=gray,fillstyle=solid](7,3){.1}
\pscircle[fillcolor=gray,fillstyle=solid](10,3){.1}

\pscircle[fillcolor=gray,fillstyle=solid](3.5,4){.1}
\pscircle[fillcolor=gray,fillstyle=solid](6,4){.1}
\pscircle[fillcolor=gray,fillstyle=solid](8,4){.1}
\pscircle[fillcolor=gray,fillstyle=solid](9.5,4){.1}
\pscircle[fillcolor=gray,fillstyle=solid](11,4){.1}

\pscircle[fillcolor=gray,fillstyle=solid](2,5){.1}
\pscircle[fillcolor=gray,fillstyle=solid](4.5,5){.1}
\pscircle[fillcolor=gray,fillstyle=solid](6.5,5){.1}
\pscircle[fillcolor=gray,fillstyle=solid](9.5,5){.1}

\pscircle[fillcolor=gray,fillstyle=solid](2,6){.1}
\pscircle[fillcolor=gray,fillstyle=solid](5,6){.1}
\pscircle[fillcolor=gray,fillstyle=solid](7.5,6){.1}
\pscircle[fillcolor=gray,fillstyle=solid](10.5,6){.1}

\psline[linecolor=gray,linewidth=1.5pt,linestyle=dashed]{->}(2,6)(2,5)
\psline[linecolor=gray,linewidth=1.5pt,linestyle=dashed]{->}(2,5)(3.5,4)
\psline[linecolor=gray,linewidth=1.5pt,linestyle=dashed]{->}(3.5,4)(2.5,3)
\psline[linecolor=gray,linewidth=1.5pt,linestyle=dashed]{->}(2.5,3)(2.5,2)
\psline[linecolor=gray,linewidth=1.5pt,linestyle=dashed]{->}(2.5,2)(2,1)
\psline[linecolor=gray,linewidth=1.5pt,linestyle=dashed]{->}(2,1)(2.5,0)

\psline[linecolor=gray,linewidth=1.5pt,linestyle=dashed]{->}(5,6)(4.5,5)
\psline[linecolor=gray,linewidth=1.5pt,linestyle=dashed]{->}(4.5,5)(3.5,4)

\psline[linecolor=gray,linewidth=1.5pt,linestyle=dashed]{->}(5,3)(5,2)
\psline[linecolor=gray,linewidth=1.5pt,linestyle=dashed]{->}(5,2)(4,1)
\psline[linecolor=gray,linewidth=1.5pt,linestyle=dashed]{->}(4,1)(2.5,0)

\psline[linecolor=gray,linewidth=1.5pt,linestyle=dashed]{->}(7.5,6)(6.5,5)
\psline[linecolor=gray,linewidth=1.5pt,linestyle=dashed]{->}(6.5,5)(6,4)
\psline[linecolor=gray,linewidth=1.5pt,linestyle=dashed]{->}(6,4)(7,3)
\psline[linecolor=gray,linewidth=1.5pt,linestyle=dashed]{->}(7,3)(7,2)
\psline[linecolor=gray,linewidth=1.5pt,linestyle=dashed]{->}(7,2)(6.5,1)
\psline[linecolor=gray,linewidth=1.5pt,linestyle=dashed]{->}(6.5,1)(7.5,0)

\psline[linecolor=gray,linewidth=1.5pt,linestyle=dashed]{->}(8,4)(7,3)
\psline[linecolor=gray,linewidth=1.5pt,linestyle=dashed]{->}(8.5,2)(8.5,1)
\psline[linecolor=gray,linewidth=1.5pt,linestyle=dashed]{->}(8.5,1)(7.5,0)

\psline[linecolor=gray,linewidth=1.5pt,linestyle=dashed]{->}(10.5,6)(9.5,5)
\psline[linecolor=gray,linewidth=1.5pt,linestyle=dashed]{->}(9.5,5)(9.5,4)
\psline[linecolor=gray,linewidth=1.5pt,linestyle=dashed]{->}(9.5,4)(10,3)
\psline[linecolor=gray,linewidth=1.5pt,linestyle=dashed]{->}(10,3)(10,2)
\psline[linecolor=gray,linewidth=1.5pt,linestyle=dashed]{->}(10,2)(10.5,1)
\psline[linecolor=gray,linewidth=1.5pt,linestyle=dashed]{->}(10.5,1)(7.5,0)
\psline[linecolor=gray,linewidth=1.5pt,linestyle=dashed]{->}(11,4)(10,3)
\rput(9,6.3){$(x,t)$}
\end{pspicture}
\caption{Howard's model and it's dual. The black points are open points, 
the gray points are the points of the dual
process and the gray (dashed) paths are the dual paths. The cluster $C(x,t)$ is enclosed 
between the left dual path and the right dual path and we have $L(x,t) = 5$.}
\label{fig:GRSDual}
\end{center}
\end{figure}

We consider $\R^2$ as space-time plane, i.e, space and time are measured along $X$ axis and $Y$ axis respectively.  The `depth' of the cluster $C(x,t)$, i.e., the time length that the cluster survived
 (see Figure \ref{fig:GRSDual}), is defined as 
$$
L(x,t) := \max\{ l \geq 0: h^l(y,s) = (x,t) \text{ for some }(y,s) \in C(x,t)\}.
$$ 
For any set $A$, we define $\# A$ to be the cardinality of $A$. It follows from Theorem \ref{thm:GRSTreeForest} that both the random variables, $\#C(x,t)$ and $L(x,t)$, are finite almost surely.

We consider the tree $T(x,t)$ obtained from the sub-graph $(C(x,t), E(x,t))$ with the natural graph metric. For $n\geq 1$,
$T_n(x,t)$ denotes the scaled tree where distances are scaled by $1/n$. 
 More formally, let $$T(x,t) := (C(x,t), E(x,t),  \mathbf{1} )$$ denote the discrete 
  weighted tree formed by the cluster $C(x,t)$,
   the edge set $E(x,t)$ and the weight function $ \mathbf{1} $, where $\mathbf{1}$ attaches costant weight $1$ to each edge in $E(x,t)$.  For $n\geq 1$
let $T_n(x,t) := (C(x,t), E(x,t),  \mathbf{1} _n)$  denote the scaled tree where 
$  \mathbf{1} _n(e) = 1/n$ for all $e \in E(x,t)$. For each $n \geq 1$, the random object $T_n(x,t)$ can be regarded as a random metric space where the distance between any two vertices (in $C(x,t)$) is given by the sum of the edge weights along the unique path between them. Note that the distribution of $T_n(x,t)$ does not depend on the vertex $(x,t) \in \Z^2$.

The main result of this paper is that, the conditional distribution of $T_n(0,0)$ given that $\{L(0,0) \geq n\}$ converges to a continuum random tree as $n \to \infty$. In order to state our result we need to describe the relevant topology briefly. 
  
Gromov-Hausdorff topology is a common way to define a topology (even a metric)
on a space of compact metric spaces. This had been introduced by Gromov \cite{M99}. 
For any two compact subsets $K$ and $K^\prime$ of a metric space $(X, d)$, the 
Hausdorff distance $d_{{\cal H}}(K, K^\prime)$ is given by:
\begin{equation*}
d_{\cal H}(K, K^\prime) := \inf\{\epsilon > 0 : K \subseteq U_{\epsilon}(K^\prime) \text{ and }K^\prime \subseteq U_{\epsilon}(K)\}
\end{equation*}
where for $\epsilon > 0$, the set $K_\epsilon$ is defined as 
$K_\epsilon := \{x \in X : d(x, K) \leq \epsilon\}$.
The Gromov-Hausdorff distance between two compact metric spaces $(X,d)$ and $(X^\prime,d^\prime)$ is defined by 
\begin{equation}
 \label{eq:GHdistMetricSpace}
 d_{\text{GH}}((X,d),(X^\prime,d^\prime))  := \inf\{d_{\cal H}(\phi(X), 
 \phi^\prime(X^\prime))\}
\end{equation}
where infimum is taken over 
all possible choices of isometric embeddings $\phi, \phi^\prime$ of the metric spaces $(X,d)$ and $(X^\prime,d^\prime)$ into a common metric space.
Let $\mathbb{M}$ be the set of all isometry equivalence classes of compact metric spaces endowed with the \textit{Gromov-Hausdorff} metric. 
 It is known that this distance turns $\mathbb{M}$
into a Polish space.
We mention here that given a metric space $(X, d)$, we write
$[X, d]$ to denote the isometry equivalence classes of $(X, d)$,
 and frequently use the notation $X$ for either $(X, d)$ or $[X, d]$ when there is no risk of ambiguity.

We show that as $n \to \infty$, the  conditional distribution of $T_n(0,0)$ given $\{L(0,0) > n\}$  converges to a continuum random tree, which we denote by $\TT$. 
Throughout this paper the notation $ \Rightarrow$ is used to denote convergence in distribution.
\begin{theorem}
\label{thm:ScalingLimCRTJt}
As $n \to \infty $, we have
$$
T_n(0,0) | \{L(0,0) \geq n\} \Rightarrow \TT,
$$
where convergence in distribution holds with respect to the Gromov-Hausdorff topology on $\mathbb{M}$.
\end{theorem}
We comment that the limiting continuum random tree $\TT$ is different from the Brownian CRT and to the best of our knowledge this limiting tree has not been studied earlier in the literature.  
In the next section, we construct $\TT$ and  explain the difference of $\TT$   from the one that  was surmised in \cite{A93}.
Slightly  related random tree has been studied  in \cite{GSW16}. In Section \ref{sec:Conclusion} we mention about the work in \cite{GSW16} in more detail. 
We don't have the complete understanding of the limiting object $\TT$ yet. Most notably 
we do not know what is the distribution of the contour walk of $\TT$ as a process (see Section \ref{sec:Conclusion} for more details). 

Scaling of conditioned tributaries are of importance because of their relations with different scaling laws empirically observed for river networks. Regarding this, we should mention that the conditional tributary for Howard's model was studied earlier and  it was proved that the joint distribution of $(L(x,t), (\# C(x,t)))$ has a regularly varying tail (Theorem 1.4 \cite{RSS16}). In order to describe another result for conditional tributary, let us describe another drainage network model.

Let $\Z^2_{\text{even}} := \{(x,t) : x, t \in \Z, x + t \text{ even}\}$
 be the oriented lattice. Each $(x,t) \in \Z^2_{\text{even}}$ acts as a source, and starting from spatial location $x$ at time $t$ water flows to location $(x+1)$ or $(x-1)$ at time $t+1$ with equal probability. 
Formally,  consider $\{ b_{(x,t)}: (x,t) \in \Z^2_{\text{even}}\}$, a collection of 
i.i.d. random variables taking  values $+1$ and $-1$ with equal probabilities. 
For $(x,t) \in \Z^2_{\text{even}}$, let $h_{\text{RW}}(x,t) := (x + b_{(x,t)}, t+1)$
and for $k \geq 1$,  $k$-th step is given by $h^k_{\text{RW}}(x,t) := h_{\text{RW}}(h^{k-1}_{\text{RW}}(x,t))$ where $h^0_{\text{RW}}(x,t) = (x,t)$. 
For $(x,t) \in \Z^2_{\text{even}}$, we observe that 
the process $\{h^{k}_{\text{RW}}(x,t)(1) : k \geq 0\}$
is a one-dimensional simple symmetric random walk starting from location $x$
at time $t$. Hence the random graph formed by the vertex set $\Z^2_{\text{even}}$ 
and the edge set $E := \{\langle (x,t),h_{\text{RW}}(x,t)\rangle:(x,t) \in \Z^2_{\text{even}}\}$ may be viewed as a graphical representation of a system of $1$-dimensional coalescing simple symmetric random walks starting from every point of  $\Z^2_{\text{even}}$. This is known as Scheidegger's model of drainage network \cite{S67}.

For $(x,t) \in \Z^2_{\text{even}}$, 
the tributary $C_{\text{RW}}(x,t)$ and its depth $L_{\text{RW}}(x,t)$ are defined similarly 
as in the case of Howard's model. Nguyen proved the following (Theorem 1 of \cite{N90}):
\begin{theorem}
\label{thm:HackNguyen}
There exists $c > 0$ such that
$$
\lim_{n \to \infty} n^{-3/2} \E(\# C_{\text{RW}}(0,0) | L_{\text{RW}}(0,0) = n)  = c.
$$
\end{theorem}
Clearly, Scheidegger's model has no cycle and in a similar way we can consider  $T^{\text{RW}}(0,0)$, which is the subtree obtained from the vertices $\{(x,t) \in \Z^2_{\text{even}} : t \geq 0\}$ and the associated edges lying in the  connected component of $(0,0)$ with the usual graph distance. For $n \geq 1$, let $T^{\text{RW}}_n(0,0)$ denote the scaled tree where distances are scaled by $1/n$.
Aldous predicted that the conditional distribution of $T^{\text{RW}}_n(0,0)$ given that $\{L_{\text{RW}}(0,0) = n\}$
converges to a CRT as $n\to \infty$ (see Subsection 4.2 Page 275 of \cite{A93}). With a modification, \cite{A93}
also provides fairly accurate description of the limiting CRT. In the same paper in Section 3, Aldous provided 
finite dimensional convergence conditions and  a tightness condition as a tool to study convergence to a continuum random tree.
Specifically for the (predicted) scaling limit of conditional tributary coming from Scheidegger's model, he further commented that:
\begin{quote}
``It is intuitively clear that rescaling the edges of $T^{\text{RW}}_n(0,0)$ by $1/n$, these random trees converge to the continuum random tree defined as above $\cdots$. This could be proved using the ideas of Nguyen \cite{N90} and Section 3." 
\end{quote}
Though the above comment suggests that finding the scaling limit is trivial in this case, but there are some important issues which we now try to highlight.

Long before, Arratia \cite{A79} observed that the system of coalescing simple symmetric random walks starting from every point on the oriented lattice has a natural dual represented  by coalescing simple symmetric random walks starting from every point on $\Z^2_{\text{odd}}:= \{(x,t) \in \Z^2: x + t \text{ odd }\}$ progressing in the backward direction of time (see Figure \ref{fig:SSRWDual}).
 To the best of our knowledge, Nguyen was the first to observe that  survival of $C^{\text{RW}}(0,0)$ can be studied through coalescence of the dual random walks starting from $(-1,0)$ and $(1,0)$.
 The dual random walks starting from $(1,0)$ and $(-1,0)$ are independent till 
 they meet. This property of independence till the time of coalescence, which is very specific
 to this model, allowed Nguyen to use the following functional limit theorem of Kaigh \cite{K76}
 to prove the result.
 \begin{figure}[htb]
\begin{center}
\scalebox{.7}{
\begin{pspicture}(0,0)(12,6)
\pscircle[fillcolor=black,fillstyle=solid](0,0){.1}
\pscircle[fillcolor=black,fillstyle=solid](2,0){.1}
\pscircle[fillcolor=black,fillstyle=solid](4,0){.1}
\pscircle[fillcolor=black,fillstyle=solid](6,0){.1}
\pscircle[fillcolor=black,fillstyle=solid](8,0){.1}
\pscircle[fillcolor=black,fillstyle=solid](10,0){.1}
\pscircle[fillcolor=black,fillstyle=solid](12,0){.1}

\pscircle[fillcolor=black,fillstyle=solid](0,2){.1}
\pscircle[fillcolor=black,fillstyle=solid](2,2){.1}
\pscircle[fillcolor=black,fillstyle=solid](4,2){.1}
\pscircle[fillcolor=black,fillstyle=solid](6,2){.1}
\pscircle[fillcolor=black,fillstyle=solid](8,2){.1}
\pscircle[fillcolor=black,fillstyle=solid](10,2){.1}
\pscircle[fillcolor=black,fillstyle=solid](12,2){.1}

\pscircle[fillcolor=black,fillstyle=solid](0,4){.1}
\pscircle[fillcolor=black,fillstyle=solid](2,4){.1}
\pscircle[fillcolor=black,fillstyle=solid](4,4){.1}
\pscircle[fillcolor=black,fillstyle=solid](6,4){.1}
\pscircle[fillcolor=black,fillstyle=solid](8,4){.1}
\pscircle[fillcolor=black,fillstyle=solid](10,4){.1}
\pscircle[fillcolor=black,fillstyle=solid](12,4){.1}

\pscircle[fillcolor=blue,fillstyle=solid](0,3){.1}
\pscircle[fillcolor=blue,fillstyle=solid](2,3){.1}
\pscircle[fillcolor=blue,fillstyle=solid](4,3){.1}
\pscircle[fillcolor=blue,fillstyle=solid](6,3){.1}
\pscircle[fillcolor=blue,fillstyle=solid](8,3){.1}
\pscircle[fillcolor=blue,fillstyle=solid](10,3){.1}
\pscircle[fillcolor=blue,fillstyle=solid](12,3){.1}

\pscircle[fillcolor=blue,fillstyle=solid](0,1){.1}
\pscircle[fillcolor=blue,fillstyle=solid](2,1){.1}
\pscircle[fillcolor=blue,fillstyle=solid](4,1){.1}
\pscircle[fillcolor=blue,fillstyle=solid](6,1){.1}
\pscircle[fillcolor=blue,fillstyle=solid](8,1){.1}
\pscircle[fillcolor=blue,fillstyle=solid](10,1){.1}
\pscircle[fillcolor=blue,fillstyle=solid](12,1){.1}

\pscircle[fillcolor=black,fillstyle=solid](1,1){.1}
\pscircle[fillcolor=black,fillstyle=solid](3,1){.1}
\pscircle[fillcolor=black,fillstyle=solid](5,1){.1}
\pscircle[fillcolor=black,fillstyle=solid](7,1){.1}
\pscircle[fillcolor=black,fillstyle=solid](9,1){.1}
\pscircle[fillcolor=black,fillstyle=solid](11,1){.1}

\pscircle[fillcolor=black,fillstyle=solid](1,3){.1}
\pscircle[fillcolor=black,fillstyle=solid](3,3){.1}
\pscircle[fillcolor=black,fillstyle=solid](5,3){.1}
\pscircle[fillcolor=black,fillstyle=solid](7,3){.1}
\pscircle[fillcolor=black,fillstyle=solid](9,3){.1}
\pscircle[fillcolor=black,fillstyle=solid](11,3){.1}

\pscircle[fillcolor=blue,fillstyle=solid](1,4){.1}
\pscircle[fillcolor=blue,fillstyle=solid](3,4){.1}
\pscircle[fillcolor=blue,fillstyle=solid](5,4){.1}
\pscircle[fillcolor=blue,fillstyle=solid](7,4){.1}
\pscircle[fillcolor=blue,fillstyle=solid](9,4){.1}
\pscircle[fillcolor=blue,fillstyle=solid](11,4){.1}

\pscircle[fillcolor=blue,fillstyle=solid](1,2){.1}
\pscircle[fillcolor=blue,fillstyle=solid](3,2){.1}
\pscircle[fillcolor=blue,fillstyle=solid](5,2){.1}
\pscircle[fillcolor=blue,fillstyle=solid](7,2){.1}
\pscircle[fillcolor=blue,fillstyle=solid](9,2){.1}
\pscircle[fillcolor=blue,fillstyle=solid](11,2){.1}

\pscircle[fillcolor=blue,fillstyle=solid](1,0){.1}
\pscircle[fillcolor=blue,fillstyle=solid](3,0){.1}
\pscircle[fillcolor=blue,fillstyle=solid](5,0){.1}
\pscircle[fillcolor=blue,fillstyle=solid](7,0){.1}
\pscircle[fillcolor=blue,fillstyle=solid](9,0){.1}
\pscircle[fillcolor=blue,fillstyle=solid](11,0){.1}

\psline[linewidth=1pt]{->}(0,0)(1,1)
\psline[linewidth=1pt]{->}(2,0)(1,1)
\psline[linewidth=1pt]{->}(4,0)(3,1)
\psline[linewidth=1pt]{->}(6,0)(7,1)
\psline[linewidth=1pt]{->}(8,0)(7,1)
\psline[linewidth=1pt]{->}(10,0)(9,1)
\psline[linewidth=1pt]{->}(12,0)(11,1)

\psline[linewidth=1pt, linecolor=blue, linestyle=dashed]{->}(1,4)(2,3)
\psline[linewidth=1pt, linecolor=blue, linestyle=dashed]{->}(2,3)(1,2)
\psline[linewidth=1pt, linecolor=blue, linestyle=dashed]{->}(1,2)(2,1)
\psline[linewidth=1pt, linecolor=blue, linestyle=dashed]{->}(2,1)(3,0)
\psline[linewidth=1pt, linecolor=blue, linestyle=dashed]{->}(3,2)(2,1)

\psline[linewidth=1pt, linecolor=blue, linestyle=dashed]{->}(3,4)(4,3)
\psline[linewidth=1pt, linecolor=blue, linestyle=dashed]{->}(4,3)(5,2)
\psline[linewidth=1pt, linecolor=blue, linestyle=dashed]{->}(5,2)(6,1)
\psline[linewidth=1pt, linecolor=blue, linestyle=dashed]{->}(6,1)(5,0)
\psline[linewidth=1pt, linecolor=blue, linestyle=dashed]{->}(4,1)(5,0)

\psline[linewidth=1pt, linecolor=blue, linestyle=dashed]{->}(3,4)(4,3)
\psline[linewidth=1pt, linecolor=blue, linestyle=dashed]{->}(4,3)(5,2)

\psline[linewidth=1pt, linecolor=blue, linestyle=dashed]{->}(5,4)(6,3)
\psline[linewidth=1pt, linecolor=blue, linestyle=dashed]{->}(6,3)(7,2)
\psline[linewidth=1pt, linecolor=blue, linestyle=dashed]{->}(7,2)(6,1)
\psline[linewidth=1pt, linecolor=blue, linestyle=dashed]{->}(8,3)(7,2)
\psline[linewidth=1pt, linecolor=blue, linestyle=dashed]{->}(7,4)(6,3)
\psline[linewidth=1pt, linecolor=blue, linestyle=dashed]{->}(9,2)(8,1)
\psline[linewidth=1pt, linecolor=blue, linestyle=dashed]{->}(8,1)(9,0)

\psline[linewidth=1pt, linecolor=blue, linestyle=dashed]{->}(9,4)(10,3)
\psline[linewidth=1pt, linecolor=blue, linestyle=dashed]{->}(10,3)(11,2)
\psline[linewidth=1pt, linecolor=blue, linestyle=dashed]{->}(11,2)(10,1)
\psline[linewidth=1pt, linecolor=blue, linestyle=dashed]{->}(10,1)(11,0)
\psline[linewidth=1pt, linecolor=blue, linestyle=dashed]{->}(9,2)(8,1)
\psline[linewidth=1pt, linecolor=blue, linestyle=dashed]{->}(11,4)(10,3)

\psline[linewidth=1pt]{->}(1,1)(0,2)
\psline[linewidth=1pt]{->}(3,1)(4,2)
\psline[linewidth=1pt]{->}(5,1)(4,2)
\psline[linewidth=1pt]{->}(7,1)(8,2)
\psline[linewidth=1pt]{->}(9,1)(10,2)
\psline[linewidth=1pt]{->}(11,1)(12,2)

\psline[linewidth=1pt]{->}(0,2)(1,3)
\psline[linewidth=1pt]{->}(2,2)(3,3)
\psline[linewidth=1pt]{->}(4,2)(3,3)
\psline[linewidth=1pt]{->}(6,2)(5,3)
\psline[linewidth=1pt]{->}(8,2)(9,3)
\psline[linewidth=1pt]{->}(10,2)(9,3)
\psline[linewidth=1pt]{->}(12,2)(11,3)

\psline[linewidth=1pt]{->}(1,3)(0,4)
\psline[linewidth=1pt]{->}(3,3)(2,4)
\psline[linewidth=1pt]{->}(5,3)(4,4)
\psline[linewidth=1pt]{->}(7,3)(8,4)
\psline[linewidth=1pt]{->}(9,3)(8,4)
\psline[linewidth=1pt]{->}(11,3)(12,4)
\end{pspicture}}
\caption{Coalescing simple symmetric random walk on $\Z^2_{\text{even}}$ and its dual. Black and blue dots represent points in $\Z^2_{\text{even}}$ and $\Z^2_{\text{odd}}$ respectively. The black edges represent
forward edges in $G$ and the blue edges represent dual edges in $\widehat{G}$. The dual graph 
represents system of coalescing simple symmetric random walks progressing in the backward 
direction of time.}
\label{fig:SSRWDual}
\end{center}
\end{figure}

\begin{theorem} 
\label{thm:Kaigh}
Random walk excursion converges to the Brownian excursion.
\end{theorem} 
The above theorem  is central in establishing Theorem \ref{thm:HackNguyen}. Now going back to the question of finding the metric space limit of $T^{\text{RW}}_n(0,0)$, in order to obtain fnite dimensional convergences, we should observe that, given $\{L^{\text{RW}}(0,0) = n\}$, the subtrees spanned by the uniformly picked vertices from the cluster $C^{\text{RW}}(0,0)$,  are  no longer given by simple random walk paths and one has to prove a suitable version of Theorem \ref{thm:Kaigh}, which is a considerable tusk.  More importantly, while dealing with models where paths are no longer independent and  have interactions, e.g., Howard's model, Nguyen's method does not help.

 Roy et al. \cite{RSS16} showed that there exists a dual process for Howard's model such that survival of the cluster can be studied through coalescence of  two neighbouring dual paths. It is important to observe that the dual process does not have the same distribution as time reversed forward process and we have  limited understanding about the distribution of the dual process. Our proof uses the fact that under diffusive scaling, Howard's model and it's dual jointly converge in distribution to the Brownian web and it's dual (Theorem 2.6 in \cite{RSS16}). Using  joint convergence and non-crossing property of paths, we show that the coalescing time of the scaled paths also converge to the coalescing time of the limiting Brownian paths. This helps us in proving metric space convergence.  In this sense the approach taken in this paper is much more robust.
 
We should mention here that Brownian web appears as universal scaling 
limit for several discrete network models, e.g., 
discrete directed spanning forest model \cite{RSS16a} 
and collection of rightmost infinite open paths starting from all percolating points for supercritical oriented percolation \cite{SS13}. For both these models, paths have  complex interactions among themselves, but the same scaling limit should hold.
 On the other hand it is important to remember that \cite{N90}
 deals with `degenerate' conditional set up $\{L(x,t) = n\}$ and we are dealing with non-degenerate conditional set up $\{L(0,0) \geq n\}$.

We would like to point out that 
this work has at least one other motivation. It has been empirically observed that river networks are structurally self similar and satisfy various scaling laws. Study of these laws and understanding the reasons behind their existence are at the core of hydrology and for this,  understanding  the behaviour of the scaled tributaries are extremely important. Ferrarri et al. \cite{FFW05}predicted the same:
\begin{quote}
``The convergence results here may lead to
rigorous/alternative verification of some of the scaling theory for those (drainage) networks."
\end{quote}
 Horton-Strahler ordering \cite{H45} is such an empirically observed scaling relation which represents scale invariance of a natural dendritic structure. Consider a binary rooted tree and assign order $1$ to each leaf. The order of an internal vertex having children with orders $i$ and $j$ 
respectively is given as 
  $$
  k = i\vee j + \mathbf{1}_{i=j} \text{ where }\mathbf{1}_A \text{ denotes indicator function for set }A.
  $$ 
Loosely speaking order of a branch denotes it's relative importance in tree hierarchy.
Sequence of connected vertices of the same order is called \textit{branch}.
 
While studying river streams, Horton \cite{H45} observed that
$N_{k+1}/N_{k} \approx R$, $3\leq R\leq 5$.
This regularity has been strongly corroborated in hydrology (\cite{SH66}, \cite{K93}, \cite{T96}, \cite{RR97})
and referred as Horton's law. 
The only rigorous result on validity of this law for drainage networks
was that of Shreve \cite{SH67}, 
who demonstrated that for a uniform distribution of rooted binary trees
 with $n$ leaves, the ratio $N_{k+1}/N_k$ converges to $4$ as $n \to \infty$. To the best of our knowledge,  there are only two other examples of binary random trees for which Horton self-similarity has been rigorously proved: the tree representation of a critical Galton 
 Watson binary branching process and tree representation of a Kingman's coalescent process.
  It is important to observe that for all these three models, there are no `space constraints' as such.  

In the context of drainage networks, any attempt to study Horton's law requires a
complete understanding of the branching structure which has complex dependencies due to space constraint. Since $T(0,0)$, the discrete tree  obtained from the river delta $C(0,0)$, contains all branching informations, it is hoped that finding this scaling limit may help in understanding Horton's law. Note that, for Howard's model $T(0,0)$ is no longer a binary tree, but Horton-Strahler ordering can be extended for a general tree in a similar way.

This paper is organized as follows. In the next section we 
construct the limiting CRT $\TT$
and explain the difference of $\TT$ than the  one described in \cite{A93}. In Section 
\ref{sec:Bweb} we introduce the Brownian web and its dual and use them to prove that $\TT$ is compact almost surely.
 In Section \ref{sec:PfThm} we describe a dual process for Howar'd model  and prove Theorem \ref{thm:ScalingLimCRTJt}. In the concluding section, i.e., Section \ref{sec:Conclusion} we make some remark about universality of our proof and present some further questions on properties of $\TT$.   

\section{Construction of the limiting CRT's}
\label{sec:TreeConstruction}

In this section we construct  the limiting CRT $\TT$.
We first observe that for $(x,t) \in \Z^2$, joining the successive 
steps $h^k(x,t),h^{k+1}(x,t) : k \geq 0$ by linear segments gives 
a continuous path $\pi^{(x,t)}$ starting from $x$ at time $t$ and
 moving in the forward direction of time as time is measured along the $Y$-axis. In other words
$\pi^{(x,t)} \in C[t,\infty)$ is such that $\pi^{(x,t)}(t+k) = h^k(x,t)(1)$
for all $k \in \N\cup\{0\}$. Here and subsequently for $(y,s) \in \R^2$, $\pi^{(y,s)}$  denotes
an element in $C[s,\infty)$ with $\pi^{(y,s)}(s) = y$. Graph distance between any two vertices in $C(x,t)$ can be interpreted as  sum of the times taken by the paths starting from each of the two vertices to coalesce. 
In the following, we use this notion to obtain a tree-like metric space 
from a given collection of paths satisfying  certain conditions
(described below in Definition \ref{def:TreeLikeMetricFrmCoalPaths}).
This notion will be used in the construction of $\TT$.
We start with a definition of  tree-like metric space or real tree.
\begin{definition}
\label{def:TreeLikeMetric}
A metric space $(X, d)$ is called   a real tree or  a $\R$-tree 
if for all $x, y \in X$ 
\begin{itemize}
 \item[(i)] there exists a unique geodesic from $x$ to $y$, i.e., there exists a unique isometry $f_{x,y}: [0, d(x,y)] \to X$ such that 
 $f_{x,y}(0) = x$ and $f_{x,y}(d(x,y)) = y$. The image of $f_{(x,y)}$ is denoted 
 by $[[x,y]]$;
 
 \item[(ii)] the only non-self-intersecting path from $x$ to $y$ is 
 $[[x,y]]$, i.e., if $q : [0,1]\to X$ is continuous and injective
  such that $q(0) = x$ and $q(1) = y$, then $
 q([0,1]) = [[xy]]$.
\end{itemize}
\end{definition} 
First condition says that $X$ is a geodesic space, and the second is a tree property that there is a unique way to travel between any two points without backtracking.
Often it is assumed that a real tree is compact. 
However for many purposes this assumption is not necessary and in
 this paper we do not make such an assumption.
  We mention here that the metric completion of an $\R$-tree is also an $\R$-tree (see \cite{I77}). Moreover it is known that the space of all isometry equivalence classes of compact real trees is closed 
in $\mathbb{M}$ (Theorem 2.1 of \cite{L05}).

Next we describe a ``tree-like" path space and how to obtain a real tree
 out of it. Let $\Pi$ be the collection of all continuous real valued paths
  moving in the forward direction of time  with all possible starting times
  (for a formal definition see Subsection \ref{subsec:Bweb}). In other words, 
 a path $\pi \in \Pi$ with starting time $\sigma_{\pi}\in \R$
is a continuous mapping $\pi :[\sigma_{\pi},\infty) \rightarrow \R $. 
Similarly $\widehat{\Pi}$ denotes the collection of all continuous real valued paths 
moving in the backward direction of time with all possible starting times.
In what follows, paths moving in the forward 
 direction of time, will be referred to as forward paths and paths moving in the 
backward direction of time, will be referred to as backward or dual paths.
The notation $\gamma(\pi_1,\pi_2)$ defined as 
$\gamma(\pi_1,\pi_2) := \inf\{ s > \sigma_{\pi_1}\vee \sigma_{\pi_2} : \pi_1(s) = \pi_2(s)\}$
denotes  the first intersection time of the  
paths $\pi_1$ and $\pi_2$ strictly after time $\sigma_{\pi_1}\vee \sigma_{\pi_2}$.
\begin{definition}[Tree-like path space]
\label{def:TreeLikeMetricFrmCoalPaths}
 $K$, a subset of $\Pi$, is said to be \textit{tree-like} if the following conditions are satisfied:
 \begin{itemize}
 \item[(i)] $\gamma(K) := \sup\{\gamma(\pi_1, \pi_2) : \pi_1, \pi_2 \in K\}$ is finite;
 \item[(ii)] for all $\pi_1, \pi_2 \in K$, we have $(\pi_1(\sigma_{\pi_1}),\sigma_{\pi_1})\neq (\pi_2(\sigma_{\pi_2}),\sigma_{\pi_2})$ and $\pi_1(s) = \pi_2(s)$ for all $s \geq \gamma(\pi_1, \pi_2)$.
 \end{itemize}
\end{definition}
These two conditions together imply that the paths in $K$ coalesce 
as soon as they intersect and hence they are non-crossing and all of them
coalesce in finite time.
The set $\{(\pi(s),s): \pi \in K, s\in [\sigma_{\pi}, \gamma(K)]\}$
consists of the images of $\pi$ in $K$ upto time $\gamma(K)$.
We assume that for any $\pi \in K$
and any $t \geq \sigma_\pi$, the path $\{\pi(s) : s \geq t\}$ with starting time $t$ is also in $K$. 

We first observe that for a tree-like $K\subset\Pi$ and for 
any $(y_1,s_1)$ in the image set, i.e., for $(y_1,s_1) \in \{(\pi(s),s): \pi \in K, s\in [\sigma_{\pi}, \gamma(K)]\}$,
there exists a unique path $\pi^{(y_1,s_1)}$ starting at $(y_1,s_1)$ in $K$. 
For any two points $(y_1,s_1)$ and $(y_2,s_2)$ in the image set,  we define the ancestor metric $d_A((y_1,s_1),(y_2,s_2))$ as,
  \begin{align}
  \label{def:AncestorMetric}
  d_A((y_1,s_1),(y_2,s_2)) :=  2\gamma(\pi^{(y_1,s_1)},\pi^{(y_2,s_2)}) - (s_1 + s_2).
  \end{align}
%
So, $d_A((y_1,s_1),(y_2,s_2))$
adds the time taken by the two paths from their start till the time they meet. 
With a slight abuse of notation, let ${\cal M}(K)$ denote the completion 
of the metric space $( \{(\pi(s),s): \pi \in K, s\in [\sigma_{\pi}, \gamma(K)]\}, d_{\cal A})$. When $\gamma(K)$ is finite, it is not difficult to see that this metric space is tree-like and 
its completion ${\cal M}(K)$ is also an $\R$-tree. 
In exactly the same way,
this notion can be extended for a collection of tree-like
backward paths $\widehat{K}\subset \widehat{\Pi}$ to obtain a complete $\R$-tree 
denoted by ${\cal M}(\widehat{K})$.


We construct our limiting object now.
Consider two independent standard Brownian motions, $B_1, B_2 $
with $B_1(0) = B_2(0) = 0$. The random times $\tau_1$ and $\tau_2$ are defined as
\begin{align*}
\tau_1 & := \inf\{ t \geq 0: \text{ for all }s_1, s_2\in [t, t+1], (B_1(s_1) - B_2(s_1))(B_1(s_2) - B_2(s_2)) \geq 0 \} \text{ and }\\
\tau_2 & := \inf\{ t > \tau_1 : B_1(t) = B_2(t)\}. 
\end{align*}
Both $\tau_1$ and $\tau_2$ are finite with $\tau_2 > \tau_1 + 1$ a.s. 
We define $B^+, B^- \in C[0,\infty]$ for $t \in [0, \tau_2 - \tau_1]$ as,
\begin{align}
\label{def:B^+B^-}
(B^+,B^-)(t)   := 
\begin{cases}
( B_1(\tau_1 + t) -  B_1(\tau_1), B_2(\tau_1 + t) -  B_1(\tau_1)) &
\text{ if }B_1(\tau_1 + 1) > B_2(\tau_1 + t)\\
( B_2(\tau_1 + t) -  B_1(\tau_1), B_1(\tau_1 + t) -  B_1(\tau_1)) &
\text{ if }B_2(\tau_1 + 1) > B_1(\tau_1 + t).
\end{cases}
\end{align}
For $t > \tau_2 - \tau_1$ we take $(B^+,B^-)(t) := ( B_1(\tau_1 + t), B_1(\tau_1 + t))$, i.e., 
the paths $B^+$ and $B^-$ coalesce at time $\tau_2 - \tau_1$.

$\Delta = \Delta(B^+, B^-)$ denotes the region in $\R^2$ 
enclosed between the paths $B^+$ and $B^-$. Formally
\begin{align*}
\Delta := \{(x,t) : t \in (0, \tau_2 - \tau_1), B^-(t) < x < B^+(t)\}.
\end{align*} 
Now we consider countable family of coalescing backward (i.e., moving in the backward 
direction of time) Brownian motions $\{\widehat{B}_m : m \in \N \}$
starting from all rational vectors in $\Delta$ such that on hitting the forward paths $B^+$ and 
$B^-$, they follow Skorohod reflection (see \cite{}). We explain this in more detail.

From the above family, consider a backward Brownian path $\widehat{B}_m$ starting from 
$(x_m, t_m) \in \Q^2 \cap \Delta$. Consider another backward Brownian 
path $\widehat{B}_{m^\prime}$ starting from $(x_m, t_m)$ independent of the forward paths $B^+$ and $B^-$. 
Then the (joint) distribution of $(B^+, B^-, \widehat{B}_m)$ is same as that of 
$(B^+, B^-, R_{B^+, B^-}(\widehat{B}_{m^\prime}))$ where 
$R_{B^+, B^-}(\widehat{B}_{m^\prime})$ is given as 
\begin{align*}
R_{B^+, B^-}(\widehat{B}_{m^\prime})(t) := 
\begin{cases}
\widehat{B}_{m^\prime}(t) - & 0 \vee \max_{t \leq s \leq t_m} (\widehat{B}_{m^\prime}(s) - B^+(t)) \\
  & - 0 \wedge \min_{t \leq s \leq t_m} (\widehat{B}_{m^\prime}(s) - B^-(t))  \qquad \text{ for }t \in [0, \tau_2 - \tau_1)\\
\widehat{B}_{m^\prime}(t) & - \widehat{B}_{m^\prime}(\tau_2 - \tau_1)   \qquad \text{ for }t \in (-\infty, 0]
\end{cases}
\end{align*}

From the work of Soucialic et. al. it follows that the finite dimensional distributions 
are consistent and hence $\{\widehat{B}_m : m \in \N\}$ exists.
By construction, we further have that the collection $\{\widehat{B}_m : m \in \N\}$
gives a non-crossing coalescing path family with $\widehat{B}_m(0) = 0$ for all $m \in \N$.      
We denote this collection of coalescing backward paths
as $\widehat{{\cal S}} = \widehat{{\cal S}}(B^+,B^-) := \{\widehat{B}_m : m \in \N\}$. 
We observe that the collection 
$\widehat{{\cal S}}$ satisfies the conditions of Definition 
\ref{def:TreeLikeMetricFrmCoalPaths} and gives a `tree-like path space'. 
The limiting continuum random tree $\TT$ is defined to be the complete metric space
 ${\cal M}(\widehat{{\cal S}})$ where completion is taken w.r.t. the ancestor metric 
 for the tree like path space $\widehat{{\cal S}}$.

\begin{remark}
\label{rem:Conj_CRT} 
\noindent {\bf $\TT$ as conjectured in \cite{A93}:}
We remark here that a slightly different continuum random tree was surmised in \cite{A93}
as the scaling limit. 
Since Aldous was interested about the scaling limit of  the
discrete tree conditioned to have time-length exactly $n$, i.e., 
$T_n(0,0)| \{L(0,0) = n\}$,
the random region $\Delta$ was enclosed by two independent backward Brownian 
motions both starting at the origin and conditioned to meet for the first time exactly at time $-1$. We are  working with the ``non-degenerate'' conditioning $\{L(0,0) \geq n\}$, hence in our case the region $\Delta = \Delta(\widehat{B}^+,\widehat{B}^-)$ is enclosed by two independent backward Brownian motions both starting at the origin and conditioned to meet before time $-1$. 
It is important to observe that for the continuum tree as described in \cite{A93}, forward coalescing Brownian paths
coalesce with the boundary of $\Delta$ as soon as they hit the boundary 
(see Page 276 of \cite{A93}). On the other hand, from the work of Soucaliuc \textit{ et. al.} \cite{STW00}, it follows that the forward coalescing Brownian paths follow Skorohod reflection
at the boundary of $\Delta(\widehat{B}^+,\widehat{B}^-)$.
This explains the difference of the limiting $\TT$  from the predicted one in \cite{A93}.
\end{remark}
 
\section{Double Brownian web and compactness of $\TT$, $\widehat{\TT}$}
\label{sec:Bweb}

In this section, we prove that $\TT$ is compact almost surely. 
Towards this  we introduce a related random object,
called the double Brownian web, i.e., the Brownian web and its dual, denoted by $(\WW, \widehat{\WW})$ and show that 
a version of ${\cal S}(\widehat{B}^+, \widehat{B}^-)$ can be embedded in $(\WW, \widehat{\WW})$.
This allows us to use properties of $(\WW, \widehat{\WW})$ to prove compactness. The proof of Theorem \ref{thm:ScalingLimCRTJt}  also uses properties of $(\WW, \widehat{\WW})$.

\subsection{Brownian web and its dual}
\label{subsec:Bweb}

The Brownian web originated in the work of Arratia (see \cite{A79}).
Later Fontes et. al. \cite{FINR04} studied the Brownian web as a random variable taking values in an appropriate Polish space. We recall the relevant details from \cite{FINR04}.

Let $\R^{2}_c$ denote the completion of the space time plane $\R^2$ with
respect to the metric
\begin{equation*}
 \rho((x_1,t_1),(x_2,t_2)) := |\tanh(t_1)-\tanh(t_2)|\vee \Bigl| \frac{\tanh(x_1)}{1+|t_1|}
 -\frac{\tanh(x_2)}{1+|t_2|} \Bigr|.
\end{equation*}
As a topological space $\R^{2}_c$ can be identified with the
continuous image of $[-\infty,\infty]^2$ under a map that identifies the line
$[-\infty,\infty]\times\{\infty\}$ with the point $(\ast,\infty)$, and the line
$[-\infty,\infty]\times\{-\infty\}$ with the point $(\ast,-\infty)$.
A path $\pi$ in $\R^{2}_c$ with starting time $\sigma_{\pi}\in [-\infty,\infty]$
is a mapping $\pi :[\sigma_{\pi},\infty]\rightarrow [-\infty,\infty] \cup \{ \ast \}$ such that
$\pi(\infty)= \ast$ and, when $\sigma_\pi = -\infty$, $\pi(-\infty)= \ast$. 
Also $t \mapsto (\pi(t),t)$ is a continuous
map from $[\sigma_{\pi},\infty]$ to $(\R^{2}_c,\rho)$.
We then define $\Pi$ to be the space of all paths in $\R^{2}_c$ with all possible starting times in $[-\infty,\infty]$.
The following metric, for $\pi_1,\pi_2\in \Pi$
\begin{equation*}
d_{\Pi} (\pi_1,\pi_2) := |\tanh(\sigma_{\pi_1})-\tanh(\sigma_{\pi_2})|\vee\sup_{t\geq
\sigma_{\pi_1}\wedge
\sigma_{\pi_2}} \Bigl|\frac{\tanh(\pi_1(t\vee\sigma_{\pi_1}))}{1+|t|}-\frac{
\tanh(\pi_2(t\vee\sigma_{\pi_2}))}{1+|t|}\Bigr|
\end{equation*}
makes $\Pi$ a complete, separable metric space. The metric $d_{\Pi}$ is slightly different 
from the original choice in \cite{FINR04} which is somewhat less natural as explained 
in \cite{SS08}. 
\begin{remark}
\label{rem:MetricConv}
Convergence in the $(\Pi, d_{\Pi})$ metric can be described as locally uniform convergence of paths as well as convergence of starting times. Therefore, for any $ \epsilon > 0$
 and $m>0$,
we can choose $ \epsilon_1 ( = g ( \epsilon, m)) > 0 $ such that for $ \pi_1, \pi_2
\in \Pi $ with $\{(\pi_i(t),t):t\in [\sigma_{\pi_i},m]\}\subseteq [-m,m]\times [-m,m]$ for $i=1,2$ and $ d_{\Pi} (  \pi_1, \pi_2 ) < \epsilon_1 $ imply that 
$|\sigma_{\pi_1}-\sigma_{\pi_2}| \vee |\pi_1 (\sigma_{\pi_1}) - \pi_2 (\sigma_{\pi_2})| < \epsilon$ and $ \sup \{| \pi_1 (t) - \pi_2 (t) | : t \in [\sigma_{\pi_1}\vee\sigma_{\pi_2},m] \} < \epsilon $.
 \qed
\end{remark}

Let ${\mathcal H}$ denote the space of compact subsets of $(\Pi,d_{\Pi})$ equipped with the Hausdorff metric $d_{{\mathcal H}}$.
As $(\Pi, d_\Pi)$ is Polish, it follows that  $({\mathcal H},d_{{\mathcal H}})$
is also Polish. Let
$B_{{\mathcal H}}$ be the Borel  $\sigma-$algebra on the metric space $({\mathcal H},d_{{\mathcal H}})$. 
The Brownian web ${\mathcal W}$ is an $({\mathcal H}, {\mathcal B}_{{\mathcal H}})$ valued random variable characterized as (Theorem 2.1 of \cite{FINR04}):
\begin{theorem}
\label{theorem:Bwebcharacterisation}
There exists an $({\mathcal H}, {\mathcal B}_{{\mathcal H}})$ valued random variable
${\mathcal W}$ such that whose distribution is uniquely determined by
the following properties:
\begin{itemize}
\item[$(a)$] for each deterministic point $\z\in \R^2$
there is a unique path $\pi^{\z}\in {\mathcal W}$  starting from $\z$ almost surely;

\item[$(b)$] for a finite set of deterministic points $\z^1,\dotsc, \z^k \in \R^2$,
the collection $(\pi^{\z^1},\dotsc,\pi^{\z^k})$ is distributed as coalescing Brownian motions starting from $\z^1,\dotsc, \z^k$;

\item[$(c)$] for any countable deterministic dense set ${\mathcal D}\subset \R^2$, 
${\mathcal W}$ is the closure of $\{\pi^{\z}: \z\in {\mathcal D} \}$ in $(\Pi, d_{\Pi})$  almost surely.
\end{itemize}
\end{theorem}
The above theorem shows that the collection is 
almost surely determined by countably many coalescing Brownian motions.

The metric space of compact sets of backward paths is defined similarly and denoted
by $(\widehat{{\mathcal H}}, d_{\widehat{{\mathcal H}}})$. Let 
$ {\mathcal B}_{\widehat{{\mathcal H}}}$ be the corresponding Borel $\sigma$ field.
$( {\mathcal W},\widehat{{\mathcal W}})$ is a $({\mathcal H}\times
\widehat{{\mathcal H}}, {\mathcal B}_{{\mathcal H}\times
\widehat{{\mathcal H}}})$ valued random variable such that ${\mathcal W}$ and $\widehat{{\mathcal W}}$ uniquely determine each other 
with $\widehat{{\mathcal W}}$ being equally distributed as $-{\mathcal W}$,
the Brownian web rotated $180^0$ about the origin. Theorem 5.8 of \cite{FINR03} and 
Theorem 2.4 of \cite{SSS16} give characterizations of 
the double Brownian web $(\WW, \widehat{\WW})$.
Before ending this subsection we list some properties of 
$(\WW, \widehat{\WW})$ which hold almost surely
and which will be used later.

%
 
\noindent {\bf Properties of $(\WW, \widehat{\WW})$ :} 
\begin{itemize}
\item[(a)] As in \cite{FINR03}, for $({\cal W}, \widehat{{\cal W}})$ and 
$(x,t) \in \R^2$,  we define 
\begin{align}
\label{def:minmout}
\begin{split}
m_{\text{in}}(x,t) := & \lim_{\epsilon \downarrow 0}\{\text{number of paths in }{\cal W}
\text{ starting at }(y,t-\epsilon) \text{ for some }y\text{ that pass }\\
&  \text{ through } (x,t) \text{ and are disjoint in the time interval }(t-\epsilon,t)\}; \\
{m}_{\text{out}}(x,t) := & \lim_{\epsilon \downarrow 0}\{\text{number of paths in }
{\cal W} \text{ starting at }(x,t)  \text{ that are } \\
& \text{ disjoint in the time interval } (t, t+\epsilon)\}.
\end{split}
\end{align}
The type of a point $ (x,t)$ is given by $ (m_{\text{in}}(x,t), m_{\text{out}}(x,t))$. 
Similarly we define $ \widehat{m}_{\text{in}}(x,t) $ and $ \widehat{m}_{\text{out}}(x,t)$
for the dual paths. 
It is known that (see Proposition 5.12, Theorem 5.13 and Theorem 5.16 of \cite{FINR03})
\begin{itemize}
\item[(i)] $m_{\text{in}}(x,t) + 1 = \widehat{m}_{\text{out}}(x,t)$ and 
$m_{\text{out}}(x,t) - 1 = \widehat{m}_{\text{in}}(x,t)$. 

\item[(ii)] Every deterministic point $ (x,t) \in \R^2$ is of type $(0,1)$. 

\item[(iii)] For any deterministic time $t$, the type of any  point on $\R \times \{t\}$ is  one of  $(0,1),(0,2)$ and $(1,1)$ in ${\cal W}$. 
\end{itemize}

\item[(b)] Let $ D \subseteq \R^2$ be a deterministic countable dense set. 
For any $ (x,t) \in D$, there exists a unique path $ \pi^{(x,t)}  \in {\cal W}$
distributed as a Brownian motion starting from $(x,t)$.

\item[(c)] Paths in $\WW$ meet and coalesce when they meet, i.e., for each $\pi, \pi^\prime \in \WW$ we have $\gamma(\pi, \pi^\prime) < \infty$ and $\pi(s) = \pi^\prime(s)$ for all $s \geq \gamma(\pi, \pi^\prime)$ (see Proposition 3.2 of \cite{SSS16}).

\item[(d)] For every $ \pi \in {\cal W} $ and $ \widehat{\pi} \in \widehat{\cal W}$  
\begin{itemize}
\item[(i)] there does not exist any $ s ,  t \in [\sigma_{\pi}, \sigma_{\widehat{\pi}}]$ with $ (\pi (s) - \widehat{\pi}(s))
(\pi (t) - \widehat{\pi}(t)) < 0 $, i.e., no  forward path of $ {\cal W} $ crosses a dual path of $ \widehat{\cal W} $;

\item[(ii)] $ \int_{\sigma_{\pi}}^{\sigma_{\widehat{\pi}}} 
\mathbf{1}_{\{ \pi(s) =  \widehat{\pi}(s) \}} ds = 0 $, i.e., 
forward paths of $ {\cal W} $ and dual paths of $ \widehat{\cal W} $ ``spend zero Lebesgue time together'' 
(see Proposition 3.2 in \cite{SS08}).
\end{itemize}

\item[(e)] For any $ t > 0, t_0 \in \R $
\begin{itemize}
\item[(i)]  the set $ \{ \pi (t_0+t) :  \pi \in {\cal W}, \sigma_{\pi} \leq t_0 \} $  is locally finite (see Proposition 4.3 in \cite{FINR04});
\item[(ii)]  for any $x \in \{ \pi (t_0+t) :  \pi \in {\cal W}, \sigma_{\pi} \leq t_0 \}$, there exist $\widehat{\pi}_1, \widehat{\pi}_2 \in \widehat{\WW}$ both starting from
the point $(x,t_0+t)$ such that $\widehat{\pi}_2(s) \neq \widehat{\pi}_1(s)$ for all $s \in [t_0,t_0 + t)$ (follows from (a));
\item[(iii)]  for any $x \in \{ \pi (t_0+t) :  \pi \in {\cal W}, \sigma_{\pi} \leq t_0 \}$, there exist $\pi_1, \pi_2 \in {\WW}$ with $\sigma_{\pi_1} < t_0 < \sigma_{\pi_2}$ such that $\pi_1(t_0 + t) = \pi_2(t_0 + t) = x$ (follows from (d) (ii), (b) and (c)).
\end{itemize}

\item[(f)] For each point $(x,t) \in \R^2$ and any two sequences $\{x^-_n : n \geq 1\}, \{ x_n^+ : n \geq 1\} \subset $ 
such that $ x_n^- \uparrow x$ and $x^+_n \downarrow x$, consider the paths $\pi^{(x^-_n,t)}, \pi^{(x^+_n,t)}
\in {\cal W}$, starting from $(x^-_n,t), (x^+_n,t)$ respectively. The limits $\lim_{n \to \infty}\pi^{(x^-_n,t)}$ and $\lim_{n \to \infty}\pi^{(x^+_n,t)}$  exist
and do not depend on the choice of the sequences $\{x^-_n : n \geq 1\}, \{ x_n^+ : n \geq 1\}$
(see Proposition 3.2 (e) of \cite{SS08}).

\item[(g)] For $ \{ \pi^n : n \geq 1 \} \subseteq {\cal W} $ and $\tilde{\pi} \in \WW$ 
with $ d_{\Pi} (\pi^n , \tilde{\pi}) 
\to 0 $, we have that $ \gamma(\pi^n, \tilde{\pi}) \to \sigma_{\tilde{\pi}} $ as $ n \to \infty$
(see Lemma 3.4 of \cite{SS08}). 

\end{itemize}

\subsection{Compactness of $\TT$}
\label{subsec:CompactCRT}

We first show that a version of ${\cal S}(\widehat{B}^+, \widehat{B}^-)$ can be embedded
in $\WW$. We need to introduce some notations.
Let $\pi^{(0,0)} \in \WW$ denote the forward Brownian path in $\WW$ 
 starting from the origin. For ease of notation we take $\pi_0 = \pi^{(0,0)}$. 
  Let $\theta_0 \leq 0$ be defined as,
 \begin{align}
 \label{def:t0}
 \begin{split}
 \theta_0 := \inf\{ &  t\geq 0 : \text{ there exists }\pi_{1}\in \WW\text{ such that }
 \sigma_{\pi_{1}} = t, \pi_{1}(t) = \pi_{0}(t)\ \\
 &\text{ and }\pi_{1}(t+1) \neq \pi_{1}(t+1)\},\\
 \theta_1 :=  & \gamma(\pi_{0}, \pi_{1}).
 \end{split}
 \end{align}
 In other words, the point $(\pi_{0}(\theta_0),\theta_0)$
has $2$ outgoing dual paths  and $\pi_{1}$ is the `newly born' 
path in $\WW$, which starts at $(\pi_{0}(\theta_0),\theta_0)$
 and do not coalesce with $\pi_0$ within the time $(\theta_0, \theta_0 + 1)$.
 $\theta_1$ denotes the coalescing time of these two paths, 
 and we have $ \theta_1 > \theta_0 + 1$ a.s.
Define 
$$
\widehat{{\cal S}}({\pi}_0, {\pi}_1) := \{\widehat{\pi} \in \widehat{\WW} 
: \sigma_{\widehat{\pi}} > \theta_0, 
(\widehat{\pi}(\sigma_{\widehat{\pi}}),\sigma_{\widehat{\pi}}) \in \Q^2, \widehat{\pi}(\theta_0)
 =  \pi_0(\theta_0) \}.
$$
$\widehat{{\cal S}}(\pi_0, \pi_1)$  precisely represents the collection of dual paths
 in $\widehat{\WW}$ starting from all rational vectors 
 in the region enclosed between $\pi_0$ and $\pi_1$.
From property (c) of $(\WW, \widehat{\WW})$, it follows that 
$\widehat{{\cal S}}(\pi_0, \pi_1)$ satisfies the conditions of Definition \ref{def:TreeLikeMetricFrmCoalPaths}.
We consider the complete metric space ${\cal M}\bigl(\widehat{{\cal S}}(\pi_0, \pi_1)
\bigr)$ with the ancestor metric and show that it has the same distribution as $\TT$.

\begin{prop}
\label{prop:EmbeddingCRTinBweb}
We have 
$$
\TT\stackrel{d}{=} {\cal M}(\widehat{{\cal S}}(\pi_0, \pi_1)).
$$
\end{prop}  
\noindent Assuming the above proposition we first prove that  $\TT$ is compact almost surely. 
 \begin{prop}
 \label{prop:CRTCompact}
 $\TT$ is compact almost surely.
 \end{prop}
 \noindent {\bf Proof :}
 Because of Proposition \ref{prop:EmbeddingCRTinBweb}, it suffices to show that 
 the metric space, ${\cal M}\bigl({\cal S}(\widehat{\pi}_0,\widehat{\pi}_1)\bigr)$ is compact almost surely.
 Being complete, it is enough to show that it is totally bounded as well. Fix $\epsilon > 0$ and set $\delta = \delta(\omega) \in (0, ((\theta_0 - \theta_1)^{-1}\wedge \epsilon)/4)$.
 For $j \geq 1$, let $t^\delta_j := \theta_1 + j \delta$ and $j_{\text{max}} := 
 \min\{j \geq 1 : t^\delta_{j+ 1} > \theta_0 \}$.
 Because of property (e) of $(\WW, \widehat{\WW})$, for any $1 \leq j \leq j_{\text{max}}-1$
  the set
 $\{\pi(t^\delta_{j+1}) : \sigma_\pi \leq t^\delta_{j}, \pi  \in {\cal S}(\widehat{\pi}_0,\widehat{\pi}_1)\}$ is finite.  Since ${\cal M}({\cal S}(\widehat{\pi}_0,\widehat{\pi}_1))$ is the completion of $\bigl (\{(\pi(s),s) : \pi \in {\cal S}(\widehat{\pi}_0,\widehat{\pi}_1), s \in [\sigma_\pi, \theta_0]\}, d_{{\cal A}}\bigr )$, the choice of $\delta$ ensures that the collection 
 $$
 \cup_{j = 1}^{j_{\text{max}}-1}\{\pi(t^\delta_{j+1}) : \sigma_\pi \leq t^\delta_{j}, \pi
 \in {\cal S}(\widehat{\pi}_0,\widehat{\pi}_1)\}\cup\{(\widehat{\pi}_0(\theta_0), \theta_0)\}
 $$ 
 forms a finite $\epsilon$-cover for ${\cal M}\bigl({\cal S}(\widehat{\pi}_0,\widehat{\pi}_1)\bigr)$.
This completes the proof.
\qed

Now we proceed to prove Proposition \ref{prop:EmbeddingCRTinBweb}. 
We need to introduce some notations. 
We define $\pi^+, \pi^- \in C[0, \infty]$ as 
\begin{equation}
\label{eq:EmbeddingCRTinBweb}
(\pi^+(s),\pi^-(s)) := 
\begin{cases}
(\pi_0( \theta_0 + s) - \pi_0( \theta_0 ),\pi_1( \theta_0 + s) - \pi_0( \theta_0 ))
 & \text{ if }\pi_0( \theta_0 + 1) > \pi_1( \theta_0 + 1)\\
 (\pi_1( \theta_0 + s) - \pi_0( \theta_0 ),\pi_0( \theta_0 + s) - \pi_0( \theta_0 ))
 & \text{ if }\pi_0( \theta_0 + 1) < \pi_1( \theta_0 + 1).
\end{cases}
\end{equation}
Since the dual paths in $\widehat{\WW}$ starting from all rational points distributed as  
coalescing backward Brownian motions and follow Skorohod reflection at the forward
 paths in $\WW$, in order to prove Proposition \ref{prop:EmbeddingCRTinBweb},
 it suffices to show that, as elements in $C[0, \infty)\times C[0, \infty)$ 
\begin{equation}
\label{eq:DeltaBoundary}
(\pi^+,\pi^-)\stackrel{d}{=} (B^+,B^-),
\end{equation}
where $(B^+,B^-)$ is as defined in (\ref{def:B^+B^-}). 

To prove Proposition \ref{prop:EmbeddingCRTinBweb} we need to deal with two issues. First of all, 
the random time $\theta_0$ is \textit{not} a stopping time. Secondly, the properties of 
$(\WW, \widehat{\WW})$ readily gives us that for each $(x,t) \in \Q^2$, there exists
unique path in $\WW$ distributed as Brownian motion starting from $(x,t)$. 
The point $(\pi_0(\theta_0), \theta_0)$ is random and to obtain the distribution of 
$(\pi^-, \pi^+)$, we have to approximate them using skeletal Brownian paths 
starting from rational points. 
Let $B_0$ and $B_{1/n}$ denote coalescing Brownian motions starting from the points 
$(0,0)$ and $(1/n, 0)$ respectively. Let $\gamma_n := \gamma(B_0, B_{1/n})$
denote the coalescing time of these two paths. 
The following proposition is the main tool for proving Proposition \ref{prop:EmbeddingCRTinBweb},
which shows that as $n \to \infty$, the conditional distribution of $(B_{1/n}, B_0)\mid \{\gamma_n
\geq 1\}$ converges to $(B^+, B^-)$. We believe that Proposition \ref{prop:ConditionalLimit_Brownian} could be of independent interest. 

\begin{prop}
\label{prop:ConditionalLimit_Brownian}
As $n \to \infty$, we have 
$$
(B_{1/n}, B_0)\mid \{\gamma_n \geq 1\} \Rightarrow (B^+, B^-).
$$
\end{prop}    

In order to prove the conditional limit theorem, we need to introduce some notations. 
Following the construction of $(B^+,B^-)$, for general 
$f_1,f_2 \in C[0,\infty)$ with $f_1(0) = f_2(0) = 0$, we define
\begin{align*}
\begin{split}
\lambda_0(f_1,f_2) := \inf\{ & t \geq 0: (f_1(s_1) - f_2(s_1))(f_1(s_2) - f_2(s_2))
> 0 \text{ for all }s_1,s_2 \in (t,t + 1)\} \text{ and }\\
\lambda_1(f_1,f_2) := \inf\{ & t > \lambda_0 : f_1(t) = f_2(t)\}.
\end{split}
\end{align*}  
For $n\geq 1$ we define
\begin{align*}
\lambda_0^n(f_1,f_2) := \inf\{ & t \geq 0: |f_1(t) - f_2(t)| = 1/n, 
(f_1(s_1) - f_2(s_1))(f_1(s_2) - f_2(s_2)) > 0 \\
& \text{ for all }s_1,s_2 \in (t,t+1)\} \text{ and }\\
\lambda_1^n(f_1,f_2) := \inf\{ & t > \lambda_0^n : f_1(t) = f_2(t)\}.
\end{align*} 
Let $A, A^n \subset C[0,\infty)\times C[0,\infty)$ be  such that 
\begin{align}
\label{def:A_An}
A := \{&(f_1,f_2) : f_i \in C[0,\infty)\text{ with } f_i(0) = 0 \text{ for }i=1,2\text{ and }\lambda_0(f_1,f_2), \lambda_1(f_1,f_2) < \infty \}\nonumber\\
A^n := \{&(f_1,f_2) : f_i \in C[0,\infty)\text{ with } f_i(0) = 0 \text{ for }i=1,2\text{ and  }, \lambda^n_0(f_1,f_2), \lambda^n_1(f_1,f_2)< \infty \}.
\end{align}
Consider the mapping $\Gamma : C[0,\infty)\times C[0,\infty) \to C[0,\infty)\times C[0,\infty)$ 
where $\Gamma(f_1,f_2) := (f_1, f_2)$ for $(f_1,f_2)\notin A$.
For $(f_1,f_2)\in A$ with $f_1(\lambda_0) = f_2(\lambda_0) = y_0 \in \R$, we 
define 
\begin{align*}
\Gamma(f_1,f_2)(s) := 
\begin{cases}
(f_1(\lambda_0 + s) - y_0, f_2(\lambda_0 + s) - y_0) & \text{ if } s \in 
[0,\lambda_1 -\lambda_0]\text{ and } f_1(\lambda_0 + 1/2) \geq f_2(\lambda_0 + 1/2)\\
(f_2(\lambda_0 + s) - y_0, f_1(\lambda_0 + s) - y_0) & \text{ if } s \in 
[0,\lambda_1 -\lambda_0]\text{ and } f_2(\lambda_0 + 1/2) > f_1(\lambda_0 + 1/2)\\
(f_1(\lambda_0 + s) - y_0, f_1(\lambda_0 + s) - y_0) & \text{ if } s \geq 
\lambda_1 -\lambda_0.
\end{cases}
\end{align*}
Similarly for $n \geq 1$, the mapping $\Gamma^n : C(-\infty, 0]\times C(-\infty, 0] \to C(-\infty, 0]\times C(-\infty, 0]$ is defined as $\Gamma^n(f_1,f_2) := (f_1, f_2)$  for $(f_1,f_2)\notin A^n$.
For $(f_1,f_2)\in A^n$ with $ f_1(\lambda^n_0) = f_2(\lambda^n_0) = y^n_0$, we  define 
\begin{align*}
\Gamma^n(f_1,f_2)(s) := 
\begin{cases}
(f_1(\lambda_0^n + s) - y^n_0, f_2(\lambda_0^n + s) - y^n_0) & \text{ if } s \in 
[0, \lambda_1^n -\lambda_0^n]\text{ and } f_1(\lambda_0^n + s) \geq f_2(\lambda_0^n + s)\\
(f_2(\lambda_0^n + s) - y^n_0, f_1(\lambda_0^n + s) - y^n_0) & \text{ if } s \in 
[0, \lambda_1^n -\lambda_0^n]\text{ and } f_1(\lambda_0^n + s) < f_2(\lambda_0^n + s)\\
(f_1(\lambda_0^n + s) - y^n_0, f_1(\lambda_0^n + s) - y^n_0) & \text{ if } s \geq 
\lambda^n_1 -\lambda^n_0.
\end{cases}
\end{align*} 

Before we proceed further, we state the following deterministic 
corollary which will be used 
 in the proof of Proposition \ref{prop:EmbeddingCRTinBweb}. For the time being, 
 we postpone the proof of Corollary \ref{cor:DeterministicContinuousMapping}.
 
\begin{corollary}
\label{cor:DeterministicContinuousMapping}
For  $(f_1,f_2) \in C(-\infty,0]\times C(-\infty,0]$ with $\lambda_1(f_1,f_2) < \lambda_0(f_1,f_2) - 1$, we have 
$$
\lim_{n \to \infty} \Gamma^n(f_1, f_2) = \Gamma(f_1, f_2)
$$
under the product metric.
\end{corollary}

Now we are ready to prove Proposition \ref{prop:ConditionalLimit_Brownian}
 and its proof is motivated from Lemma 3.1 of \cite{B76}.

\noindent {\bf Proof of Proposition \ref{prop:ConditionalLimit_Brownian}:}
Let $B_1, B_2$ are two independent Brownian motions starting from 
the origin. We observe that almost surely, 
as pair of paths $(B_1, B_2) \in A_n$ for all $n \geq 1$.
We first show that for all $n \geq 1$, we have
\begin{align}
\label{eq:Cond_Brown_Dist_equality}
(B_{1/n},B_0)|\{ \gamma_n < - 1\} \stackrel{d}{=} \Gamma^n(B_1, B_2).
\end{align}

For $t > 0$ and $n \geq 1$, we define the event $E^n_t$ as 
\begin{align*}
\bigcup_{ s \in [0,t)\cap \Q} \Bigl \{ (B_1 - B_2)(s) >  1/n \text{ and we have:}
|(B_1 - B_2)(s^\prime)| > 0 \text{ for all } s^\prime \in [s, t\wedge(s+1)] \Bigr\}.
\end{align*}
We note that occurrence of the event $E^n_t$ implies that there exists some $s_0 <t$ 
such that $|(B_1 - B_2)(s_0)| = 1/n$ with $|(B_1 - B_2)(s^\prime)| > 0$
for all $s^\prime \in (s_0 , t\wedge(s+1))$. 
By definition we have that for all $t > 0$, the event $E^n_t$ is 
${\cal F}_t := \sigma(\{B_1(s) , B_2(s) : 0 \leq s \leq t\})$ measurable.

For any $t> 0$ and for $n\geq 1$  we have 
\begin{align}
\label{eq:EventEquality1}
\{\lambda^n_0(B_1, B_2) = t\} = &
\{|B_1(s_1) - B_2(s_1)| > 0 \text{ for all }s_1 \in [t, t+1]\}\nonumber \\
& \cap \{|B_1(t) - B_2(t)| = 1/n\}\cap(E^n_t)^c,
\end{align}
as non-occurrence of the event $E^n_t$ ensures that the random time $\lambda^n_0$ can not occur earlier than $t$ and the event  $\{|B_1(s_1) - B_2(s_1)| > 0 \text{ for all }s_1 \in [t, t+1]\}
 \cap \{ |B_1(t) - B_2(t)| = 1/n\}$ confirms it can not occur later than $t$.

Set $B_1(\lambda^n_0) \wedge B_2(\lambda^n_0) = y^n_0$. 
Using $B_1, B_2$ now we define $B^\prime_1, B^\prime_2 \in C[0,\infty)$ as 
\begin{align*}
(B^\prime_1(s), B^\prime_2(s)) := 
\begin{cases}
(B_1(\lambda^n_0 + s) - y^n_0, B_2(\lambda^n_0 + s) - y^n_0) & 
\text{ if }s \in [0, \lambda^n_1 - \lambda^n_0]\\
(B_1(\lambda^n_0 + s) - y^n_0,B_1(\lambda^n_0 + s) - y^n_0) & \text{ if }s > \lambda^n_1 - \lambda^n_0.
\end{cases}
\end{align*} 
For any Borel set $B \subseteq C[0, \infty) \times C[0, \infty)$, we obtain
\begin{align}
\label{eq:Calculation}
& \P \bigl(\Gamma^n(B_1, B_2) \in B \bigr)\nonumber\\
= & \E \Bigl( \P \bigl(\Gamma^n(B_1, B_2) \in B  \bigl | \lambda^n_0 = t \bigr ) \Bigr)\nonumber\\
= & \E \Bigl( \P \bigl(\{(B^\prime_1, B^\prime_2)(t+s): s \geq 0\} \in B  \bigl | 
 |(B_1 - B_2)(t)| = 1/n, |(B_1 - B_2)(s_1)| > 0 \\
 & \qquad \qquad \qquad \qquad \text{ for all }s_1 \in [t, t+1],
  (E^n_t)^c \bigr )\Bigr). 
\end{align}
The event $(E^n_t)$ is measurable w.r.t. the $\sigma$-field
${\cal F}_t := \sigma(\{B_1(s) , B_2(s) : 0 \leq s \leq t\})$ 
and from Markov property of $(B_1, B_2)$ 
it follows that given $\{|(B_1 - B_2)(t)| = 1/n, |(B_1 - B_2)(s_1)| > 0 
\text{ for all }s_1 \in [t, t+1], (E^n_t)^c \}$, the distribution of the 
process $\{(B^\prime_1(t+s), B^\prime_2(t+s)) : s \geq 0\}$ does not depend on $E^n_t$ and 
has the same distribution as two independent coalescing
Brownian motions starting from the points $(1/n,0)$ and $(0,0)$ respectively and conditioned 
\textit{not} to meet before time $1$. Hence from (\ref{eq:Calculation}), we obtain 
\begin{align*}
= & \E \Bigl( \P \bigl(\{(B^\prime_1, B^\prime_2)(t+s): s \geq 0\} \in B  \bigl | 
 |(B_1 - B_2)(t)| = 1/n, |(B_1 - B_2)(s_1)| > 0 \\
 & \qquad \qquad \qquad \qquad \text{ for all }s_1 \in [t, t+1],
  (E^n_t)^c \bigr )\Bigr)\\
= & \E \Bigl(\P \bigl((B^{1/n},B_0) \in B \bigl| \gamma_n >  1 \bigr)\Bigr)\\
= & \P \bigl((B^{1/n},B_0) \in B \bigl| \gamma_n >  1 \bigr).
\end{align*}
This completes the proof of (\ref{eq:Cond_Brown_Dist_equality}). 

Now fix any bounded continuous function $f : C[0,\infty) \times C[0,\infty) \mapsto \R$. 
\begin{align*}
& \lim_{n \to \infty}  \E( f(B_{1/n}, B_0) \mid \gamma_n > 1)\\
& = \lim_{n \to \infty}  \E( f(\Gamma_n(B_1,B_2 ))) \text{ by }(\ref{eq:Cond_Brown_Dist_equality})\\
& =    \E( f(\Gamma(B_1,B_2 ))) \text{ by bounded convergence theorem and by
 Corollary }\ref{cor:DeterministicContinuousMapping}\\
& =    \E( f(B^+,B^- )).
\end{align*}
Since $f$ is chosen arbitrarily, this completes the proof. 
\qed

Now we use Proposition \ref{prop:ConditionalLimit_Brownian} to prove Proposition \ref{prop:EmbeddingCRTinBweb}. 

\noindent {\bf Proof of Proposition \ref{prop:EmbeddingCRTinBweb}:}
As  observed earlier it suffices to prove (\ref{eq:DeltaBoundary}).
 Recall that  $\widehat{\pi}_0$ is the dual path in $\widehat{\WW}$ starting from the point $(0,0)$.  For $n \geq 1$ let $l_n$ and $r_n$ be rationals given by $l_n := (\lfloor n \widehat{\pi}_0(\theta_0) \rfloor)/n$ and 
 $r_n := l_n + 1/n$ where $\theta_0$ is as defined as in (\ref{def:t0}).
 We choose $\theta^n_0 \in (\theta_0 - 1/n, \theta_0)\cap \Q$ such that 
 \begin{itemize}
 \item[(i)] $\widehat{\pi}_0(s),\widehat{\pi}_1(s) \in [l_n,r_n)$ for all $s \in [\theta^n_0, \theta_0]$ where $\widehat{\pi}_1$ is the dual path in $\widehat{\WW}$ starting from $(\widehat{\pi}_0(\theta_0),\theta_0)$, 
 \item[(ii)] $\theta^n_0 > \theta_1 + 1$.
 \end{itemize}
 Since $\theta_0 > \theta_1 + 1$ almost surely and both the backward continuous paths $\widehat{\pi}_0$ and $\widehat{\pi}_1$ pass through the point $(\widehat{\pi}_0(\theta_0),\theta_0)$, such $\theta^n_0$ always exists. Clearly $\theta^n_0 \to  \theta_0 $ as $n \to \infty$.
 For all $n \geq 1$ both $(l_n,\theta^n_0)$ and $(r_n,\theta^n_0)$ are in $\Q^2$ and 
 let $\widehat{\pi}^{(l_n,\theta^n_0)}$ and $\widehat{\pi}^{(r_n,\theta^n_0)}$ denote the dual 
 paths in $\widehat{\WW}$ starting from the points $(l_n,\theta^n_0)$ and $(r_n,\theta^n_0)$
 respectively. Since the dual paths in $\widehat{\WW}$ are non-crossing, condition (ii) implies that $\widehat{\pi}^{(r_n,\theta^n_0)}(\theta^n_0-1) > \widehat{\pi}^{(l_n,\theta^n_0)}(\theta^n_0-1)$.   The backward paths $\widehat{\pi}^+_n, \widehat{\pi}^-_n \in C(-\infty, 0]$ are defined as  
 \begin{align*}
 \widehat{\pi}^+_n(-s) := \widehat{\pi}^{(r_n,\theta^n_0)}( \theta^n_0 - s ) - l_n\text{ and }
 \widehat{\pi}^-_n(-s) := \widehat{\pi}^{(l_n,\theta^n_0)}( \theta^n_0 - s ) - l_n \text{ for }s\geq 0.
 \end{align*}
  By construction we have $\widehat{\pi}^+_n(0) - \widehat{\pi}^-_n(0) = 1/n$ and 
    $\widehat{\pi}^+_n(s) > \widehat{\pi}^-_n(s)$
 for all $s \in [-1,0]$ almost surely.

Since $\widehat{\WW}$ is compact, both the sequences $\{\widehat{\pi}^{(r_n,\theta^n_0)} : n \in \N\}$
and $\{\widehat{\pi}^{(l_n,\theta^n_0)} : n \in \N\}$ must have convergent subsequences. 
As convergence of paths implies convergence of starting times as well, subsequential limits of each of these two sequences must be backward paths in $\widehat{\WW}$ 
starting at $(\widehat{\pi}_0(\theta_0),\theta_0)$.
 From the properties of $(\WW,\widehat{\WW})$ we have that there are exactly two dual paths starting from the point $(\widehat{\pi}_0(\theta_0), \theta_0)$. Hence it follows that 
\begin{align}
\label{eq:limLnRnPath}
\lim_{n\to \infty} (\widehat{\pi}^{(l_n,\theta^n_0)},\widehat{\pi}^{(r_n,\theta^n_0)})
=\begin{cases}
(\{\widehat{\pi}_0(s) : s \leq \theta_0\}, \widehat{\pi}_1)&\text{ if }\widehat{\pi}_0(\theta_0 - 1) < \widehat{\pi}_1(\theta_0 - 1)\\
(\widehat{\pi}_1, \{\widehat{\pi}_0(s) : s \leq \theta_0\})&\text{ otherwise},
\end{cases}
\end{align} 
where the limit is taken in $(\widehat{\Pi}, d_{\widehat{\Pi}})$.
This shows that the sequences of  
 dual paths, $\{\widehat{\pi}^+_n : n \in \N\}$ and $\{\widehat{\pi}^-_n : n \in \N\}$,
 almost surely converge to the dual paths $\widehat{\pi}^+$ and $\widehat{\pi}^-$ respectively.

\noindent Fix any bounded continuous function $g : C(-\infty,0]\times C(-\infty,0] \mapsto \R$.
We claim that
 \begin{equation}
 \label{eq:limPin}
 \lim_{n \to \infty}\E(g(\widehat{\pi}^+_n, \widehat{\pi}^-_n)) =
\lim_{n \to \infty}\E(g(\widehat{W}^{(1/n,0)}, \widehat{W}^{(0,0)})|\widehat{\gamma}_n < -1).
 \end{equation}
 The proof is similar to the proof of Lemma \ref{lem:ConditionalStoppingTime}.
 For completeness we present here full details.
 For each $n \geq 1$, let $O_n$ denote the event that the set of the points $\{\widehat{\pi}_0(\theta_0 - 1), \widehat{\pi}_1(\theta_0 - 1)\}$ equals the set $\{\widehat{\pi}^{(l_n,\theta^n_0)}(\theta_0 - 1), \widehat{\pi}^{(r_n,\theta^n_0)}(\theta_0 - 1)\}$.  
 From (\ref{eq:limLnRnPath}) and from property (g) of $(\WW, \widehat{\WW})$, we have that 
 $\P(O_n) \to 1$ as $n \to \infty$.
 Hence we have 
\begin{align}
\label{eq:LimEgOn}
\lim_{n \to \infty}\E(g(\widehat{\pi}^+_n, \widehat{\pi}^-_n)) = 
 \lim_{n \to \infty}\E(g(\widehat{\pi}^+_n, \widehat{\pi}^-_n)\mathbf{1}_{O_n}),
 \end{align}
 where $\mathbf{1}_{O_n}$ denotes indicator  function of the event $O_n$.
 
\noindent We need to define some more events.
For $t > 0 $ and $n \geq 1$, let $F_t$ denote the event that 
\begin{align*}
\bigcup_{ s \in (0,t)\cap \Q}\quad \bigcap_{ m \geq 1}\quad \bigl \{ & \text{there exists }
\widehat{\pi}^\prime_m \in \widehat{\WW} \text{ with }\sigma_{\widehat{\pi}^\prime_m} \geq -s_m
\text{ for some }s_m \in (0,s)\cap\Q \text{ such that }\\ 
&|\widehat{\pi}_0(-s_m) - \widehat{\pi}^\prime_m(-s_m)| < 1/m
\text{ and }\widehat{\pi}_0((-s -1 )\vee(-t)) \neq \widehat{\pi}^\prime_m((-s -1 )\vee(-t)) \bigr\}.
\end{align*}
 As $\widehat{\WW}$ is compact, the sequence of dual paths $\{\widehat{\pi}^\prime_m : m \in \N\}$ considered in the event $F_t$ must have  a subsequential limit.
Hence occurrence of the event $F_t$ ensures existence of a dual path $\widehat{\pi}_1\in \widehat{\WW}$ starting from $(\widehat{\pi}_0(-s_0),-s_0)$ for some $s_0 < t$ such that
$\widehat{\pi}_0((-s_0 -1 )\vee(-t)) \neq \widehat{\pi}_1((-s_0 -1 )\vee(-t))$.

\noindent Next for $t>0, i \in \Z$ and $s \in Q$, we define the following events :
\begin{align*}
E_t(i,n) &:= \{\widehat{\pi}_0(-t) \in [i/n,(i+1)/n)\};\\
E_t(i,n,s) &:= \{\widehat{\pi}^{(\frac{i+1}{n},s)}(s -1) > \widehat{\pi}^{(\frac{i}{n},s)}(s -1)\},
\end{align*}
 where $\widehat{\pi}^{(\frac{i+1}{n},s)}$ and $\widehat{\pi}^{(\frac{i}{n},s)}$ are the dual paths in $\widehat{\WW}$ starting from the points $(\frac{i+1}{n},s)$ and $(\frac{i}{n},s)$ respectively.
For $t > 0$ with a slight abuse of notation let ${\cal G}_t$ denote the $\sigma$-field given by $${\cal G}_t := \sigma \bigl(\{\widehat{\pi}(s) : \widehat{\pi} \in \widehat{\WW}\text{ with } \sigma_{\widehat{\pi}}\geq -t\text{ and } s \in [-t, \sigma_{\widehat{\pi}}]\}\bigr).$$
By definition for all $t > 0$, both the events $F_t$ and $E_t(i,n)$ are 
${\cal G}_t $ measurable. 
Hence from non-crossing nature of paths in $\widehat{\WW}$ it follows that for all $t > 0$ and for $t^n_0 \in (t,t+1/n)\cap \Q$, we have the following equality of  events
\begin{align}
\label{eq:EventEquality2}
& \{\theta_0 = -t\}\cap \{\theta^n_0 = -t^n_0\}\cap E_t(i,n)\cap O_n = \nonumber\\
& \qquad \qquad \qquad (F_t)^c 
\cap\{\widehat{\pi}_0(s),\widehat{\pi}_1(s) \in [i/n,(i+1)/n)\text{ for all }s\in [-t^n_0,-t]\}\cap E_t(i,n,-t^n_0)\cap O_n.
\end{align}

For $t > 0$ and for $t^n_0 \in (t,t+1/n)\cap \Q$ using (\ref{eq:LimEgOn}) and (\ref{eq:EventEquality2}) we obtain
\begin{align}
\label{eq:CalculationEmbeddingCRT1}
& \lim_{n \to \infty}\E(g(\widehat{\pi}^+_n, \widehat{\pi}^-_n)) \nonumber\\
& = \lim_{n \to \infty} \E\bigl( \E(g(\widehat{\pi}^+_n, \widehat{\pi}^-_n) \mathbf{1}_{O_n}|  \theta_0 = - t,\theta^n_0 = -t^n_0) \bigr)\nonumber\\
& = \E\bigl(\lim_{n \to \infty}  \E(g(\widehat{\pi}^+_n, \widehat{\pi}^-_n) \mathbf{1}_{O_n}|\theta_0 = - t,\theta^n_0 = -t^n_0) \bigr)\nonumber\\
&  = \E\Bigl(\lim_{n\to\infty}\sum_{i=-\infty}^{\infty} 
\E \bigl(g(\widehat{\pi}^+_n, \widehat{\pi}^-_n)\mathbf{1}_{O_n} 
\bigl| \theta_0 = - t, \theta^n_0 = -t^n_0, E_t(i,n)\bigr)\P(E_t(i,n))\Bigr)\nonumber\\
&  = \E\Bigl(\sum_{i=-\infty}^{\infty} 
\lim_{n\to\infty}\E \bigl(g(\widehat{\pi}^+_n, \widehat{\pi}^-_n)\mathbf{1}_{O_n} 
\bigl| \theta_0 = - t, \theta^n_0 = -t^n_0,E_t(i,n)\bigr)\P(E_t(i,n))\Bigr)\nonumber\\
&  = \E\Bigl(\sum_{i=-\infty}^{\infty} \lim_{n\to\infty} \E  \bigl(g(\{(\widehat{\pi}^{(\frac{i+1}{n},-t^n_0)}(-t^n_0 - s) - i/n,\widehat{\pi}^{(\frac{i}{n},-t^n_0)}(-t^n_0 -s) - i/n) : s \geq 0\}) \bigl| \nonumber\\
& \qquad   (F_t)^c, E_t(i,n,-t^n_0)\text{ and } \widehat{\pi}_0(s),\widehat{\pi}_1(s) \in [i/n,(i+1)/n)\text{ for all }s\in [-t^n_0,-t] \bigr)\P( E_t(i,n))\Bigr).
\end{align}
Next we observe that conditioned on the event that $(F_t)^c\cap  E_t(i,n,-t^n_0)\cap \{ \widehat{\pi}_0(s),\widehat{\pi}_1(s) \in [i/n,(i+1)/n)\text{ for all }s\in [-t^n_0,-t]\}$, the process $\{(\widehat{\pi}^{(\frac{i+1}{n},-t^n_0)}(-t^n_0 - s) - i/n,\widehat{\pi}^{(\frac{i}{n},-t^n_0)}(-t^n_0 -s) - i/n) : s \geq 0\}$ is independent of the $\sigma$-field  ${\cal G}_{t^n_0}$ and has the same distribution as $(\widehat{W}^{(1/n,0)},\widehat{W}^{(0,0)})|\{\widehat{\gamma}_n < -1\}$.  Since $\{\widehat{\pi}_0(s),\widehat{\pi}_1(s) \in [i/n,(i+1)/n)\text{ for all }s\in [-t^n_0,-t]\}\cap (F_t)^c$
 is ${\cal G}_{t^n_0}$ measurable, from (\ref{eq:CalculationEmbeddingCRT1}) we obtain 
\begin{align*}
&  = \E\Bigl(\sum_{i=-\infty}^{\infty} \lim_{n\to\infty} \E  \bigl(g(\{(\widehat{\pi}^{(\frac{i+1}{n},-t^n_0)}(-t^n_0 - s) - i/n,\widehat{\pi}^{(\frac{i}{n},-t^n_0)}(-t^n_0 -s) - i/n) : s \geq 0\}) \bigl| \nonumber\\
& \qquad \qquad E_t(i,n,-t^n_0) \bigr)\P(E_t(i,n))\Bigr)\\
& = \E\bigl(\sum_{i=-\infty}^{\infty} \E(g(\widehat{W}^{(1/n,0)}, \widehat{W}^{(0,0)}) \bigl|  \widehat{\gamma}_n < - 1)\P(E_t(i,n))\bigr)\nonumber\\
& = \E(g(\widehat{W}^{(1/n,0)}, \widehat{W}^{(0,0)}) |  \widehat{\gamma}_n < -1).
\end{align*}
This completes the proof of (\ref{eq:limPin}).  
Since $(\widehat{\pi}^+_n,\widehat{\pi}^-_n) \to (\widehat{\pi}^+,\widehat{\pi}^-)$ almost surely as $n\to \infty$, from Lemma \ref{lem:CRTConditionalLimit} and Corollary \ref{cor:DeterministicContinuousMapping} we have 
\begin{align*}
\E(g(\widehat{\pi}^+,\widehat{\pi}^-)) & = \lim_{n \to \infty}\E(g(\widehat{\pi}^+_n,\widehat{\pi}^-_n))\\
& =\lim_{n \to \infty}\E(g(\widehat{W}^{(1/n,0)},\widehat{W}^{(0,0)})|\widehat{\gamma}_n < -1)\\
& =\lim_{n \to \infty}\E(g(\Gamma_n(\widehat{B}_1,\widehat{B}_2)))\\
& =\E(g(\Gamma(\widehat{B}_1,\widehat{B}_2)))=\E(g(\widehat{B}^+,\widehat{B}^-)). 
\end{align*}
This completes the proof.
 \qed 

Finally we need to complete the proof of Corollary \ref{cor:DeterministicContinuousMapping}. 

\noindent {\bf Proof of Corollary \ref{cor:DeterministicContinuousMapping}:} We observe that  for all large $n$, both $\lambda^n_1(f_1,f_2)$ and $\lambda^n_0(f_1,f_2)$ are finite and $\lambda^n_1(f_1,f_2) < \lambda^n_0(f_1,f_2) - 1$. Fix $\epsilon > 0$. For ease of notations, we take $\lambda_0 = \lambda_0(f_1,f_2), \lambda_1 = \lambda_1(f_1,f_2), \lambda^n_0 = \lambda^n_0(f_1,f_2)$ and $\lambda^n_1 = \lambda^n_1(f_1,f_2)$. Since $f_1, f_2$ are both continuous with $f_1(0) = f_2(0) = 0$, we have $\lambda^n_0 \leq \lambda_0$ and $\lambda^n_1 \leq \lambda_1$ for all large $n$. Set $n_0$ such that $\lambda^n_0 \leq \lambda_0,  \lambda^n_1 = \lambda_1$ and $1/n < \epsilon/2$ for all $n \geq n_0$. Choose $\delta > 0$ such that
the following conditions hold :
\begin{itemize}
\item[(i)] $\min\{|f_1(s) - f_2(s) | : s \in [\lambda_0 - 1 - \delta, \lambda_0 - \delta] \} = \beta > 0$,
\item[(ii)] $\sup\{|f_1(x) - f_1(y)|\vee |f_2(x) - f_2(y)| : x,y \in [ \lambda_1, 0],
 |x-y| < \delta \} < \epsilon/2 $ and
\item[(iii)] $\sup\{|f_1(x) - f_2(y)|: x,y \in [\lambda_1-\delta, \lambda_1 + \delta]\} < \epsilon/2$. 
\end{itemize} 
Since $f_1$  and $f_2$ are both continuous with $f_1(\lambda_1) = f_2(\lambda_1)$ and $\lambda_1$ is strictly smaller than $ \lambda_0 - 1$, $\delta> 0$ always exists. 
Condition (ii) ensures that $\lambda^n_0 \geq \lambda_0 -  $
Because of condition (i), for all $ n > 1/\beta$ we have 
$\lambda_1^n \in (\lambda_0 - \delta, \lambda_0) $ and consequently $\lambda_1^n = \lambda_1$. Further  the conditions (ii) and (iii) mentioned above imply that 
\begin{align}
\label{CalculationGamma}
& \sup\{|f_1(\lambda_0^n - s) - f_1(\lambda_0 - s)|\vee |f_2(\lambda_0^n - s) - f_2(\lambda_0 - s)| : s \in [0, \lambda_0 - \lambda_1]\} < \epsilon 
\text{ and }\nonumber \\
& \Gamma_n(f_1, f_2)(-s) = \Gamma(f_1, f_2)(-s)\text{ for all }s\geq \lambda_0 - \lambda_1.
\end{align}
Since $\epsilon > 0$ is chosen 
arbitrarily, (\ref{CalculationGamma}) completes the proof.
\qed

The embedding discussed in Proposition \ref{prop:EmbeddingCRTinBweb} proves the following corollary:
\begin{corollary}
\label{cor:CRTRootType}
For any $\epsilon \in (0,1)$ and for $\pi_1, \pi_2 \in {\cal S}(\widehat{B}^+, \widehat{B}^-)$ with $\sigma_{\pi_1}\vee \sigma_{\pi_2} \leq - \epsilon$, almost surely there exists $s_0 = s_0(\omega) \in (-\epsilon, 0)$ such that $\pi_1(s) =  \pi_2(s)$ for all $s \geq s_0$.
\end{corollary} 
\noindent \textbf{Proof :} The proof follows from Proposition \ref{prop:EmbeddingCRTinBweb}
and from property (a) of $(\WW, \widehat{\WW})$ which ensures that $m_{\text{in}}(\widehat{\pi}_0(\theta_0), \theta_0) = 1$ almost surely.
\qed

The same argument as in Proposition \ref{prop:EmbeddingCRTinBweb} gives us the following 
corollary:
\begin{corollary} 
\label{}
  
\end{corollary}

We end this section with the following two remarks:
\begin{remark}
\label{rem:Condi_DualPaths_Dist}
Conditional that the point $(x,t)$ is a $(1,2)$ type point, i.e., one incoming path and two outgoing 
paths (see Section \ref{sec:Bweb} for a detailed definition) such that these two outgoing paths,
 denoted by $\pi^{(x,t)}_+$ and $\pi^{(x,t)}_-$, don't coalesce before time $t+1$. 
We define $(\pi^+, \pi^-) \in C[0,\infty)\times C[0,\infty)$ such that 
\begin{align*}
(\pi^+, \pi^-)(s) :=
\begin{cases}
(\pi^{(x,t)}_+(t+s) - x, \pi^{(x,t)}_-(t+s) - x) & \text{ for } 0 \leq s \leq 
\gamma(\pi^{(x,t)}_+, \pi^{(x,t)}_-) - t \\
(\pi^{(x,t)}_+(t+s) - x, \pi^{(x,t)}_+(t+s) - x) & \text{ for } s \geq 
\gamma(\pi^{(x,t)}_+, \pi^{(x,t)}_-) - t .
\end{cases}
\end{align*}
Same argument as in Proposition \ref{prop:EmbeddingCRTinBweb} gives us the distributional equality 
\begin{align}
\label{eq:BrownianMeanderDist_equality}
(\pi^+, \pi^-) \stackrel{d}{=}(B^+,  B^-).
\end{align}  
\end{remark}

\begin{remark}
\label{rem:CRTRootType}
Regarding the construction of $\TT$, it is useful to mention 
here that instead of taking closure with respect to the ancestor metric, 
one can take closure of $\widehat{{\cal S}}(B^+, B^-)$ in $(\widehat{\Pi}, d_{\widehat{\Pi}})$. 
The embedding explained in Proposition \ref{prop:EmbeddingCRTinBweb}
shows that in that case the resulting collection of dual paths will no longer satisfy
the condition of Definition \ref{def:TreeLikeMetricFrmCoalPaths}, 
as almost surely there will be (random) points with 
multiple outgoing paths. Nevertheless, the resulting object almost surely will be an  
$\R$-graph (for a formal definition and basic properties of an $\R$-graph, see \cite{BBGM15}).
For constructing the limiting real tree, one could have also followed  
\textrm{cutting} operation of this $\R$-graph as in \cite{BBGM15}.
Finally we comment that the distribution of $(\TT, \widehat{\TT})$
 does not depend on the choice of the deterministic countable dense set of $\R^2$. 
\end{remark}
 
\section{ Proof of Theorem \ref{thm:ScalingLimCRTJt}}
\label{sec:PfThm}

In this section we prove Theorem \ref{thm:ScalingLimCRTJt}.
Before embarking on the details, we sketch the idea of the proof.
Coletti  et. al. \cite{CFD09} showed that the scaled Howard's model 
observed as collection of paths, converges in distribution to the Brownian web 
(see Theorem  in \cite{CFD09}). Later, Roy et. al. constructed a dual 
network for the Howard's model and showed that under diffusive scaling, Howard's 
model and it's dual jointly converges in distribution to the Brownian web and it's 
dual (Theorem of \cite{RSS16}). We use this joint convergence in path-space topology
 and prove that a scaled pair of Howard paths together with it's coalescing time
(which is finite due to due to Theorem \ref{thm:GRSTreeForest}) jointly 
converges in distribution to a pair of coalescing Brownian paths and their coalescing time.
 This enables us to obtain convergence of metric spaces in the Gromov-Hausdorff topology.

We need to introduce some notation. 
For $(x,t) \in \Z^2$, taking the edges
$ \langle h^{k-1}(x,t),h^{k}(x,t) \rangle $  
to be straight line segments for $k \geq 1$, we parametrize the path formed by these edges as
the piece wise linear function $\pi^{(x,t)} : [t, \infty) \to \mathbb R$
such that ${\pi}^{(x,t)}(t+k) = h^{k}(x,t)(1)$ for every integer $k \geq 0$ and
 and linear in between.
Let ${\cal X} := \{\pi^{(x,t)}:(x,t) \in \Z^2_{\text{even}}\}$ denote the collection of all (forward) paths obtained from $G$.
For each $n \geq 1$ and for $\pi \in \Pi$,  the scaled path
$ \pi_n : [\sigma_{\pi}/n, \infty] \to
[-\infty, \infty]$ is given by $\pi_n (t) :=
\pi(n t)/\sqrt{n} $. 
For each $ n \geq 1 $, let ${\cal X}_n :=
\{\pi_n^{(x,t)} :(x,t) \in \Z^2_{\text{even}}\}$ be the collection of the scaled paths and 
$\bar{\cal X}_n$ is the closure of
${\cal X}_n$ in $(\Pi,d_{\Pi})$.
As we commented earlier, in order to obtain the scaling limit one has to deal 
with the complex dependencies between multiple paths. Coletti et. al. showed that 
under diffusive scaline, the Howard's model converges in distribution to the 
Brownian web. 
\begin{theorem}[Theorem 6.1 of \cite{CFD09}]
\label{th:SSRW_BW}
As $n\rightarrow \infty$,  $\bar{{\cal X}}_n$
 converges in distribution to the Brownian web ${\cal W}$ as 
 ${\cal H},{\cal B}_{{\cal H}})$ valued random variables.
\end{theorem}

Later Roy et. al. constructed a dual graph for Howard's model 
and extended the above result to obtain joint convergence to the double Brownian web 
$(\WW, \widehat{\WW})$ for Howard's model and it's dual under diffusive scaling. 
Similarly we construct a scaled family of backward/dual paths obtained from $\widehat{G}$.

We need to explain the dual graph first. 
For $(x,t)\in \Z^2_{\text{odd}}$,
the dual path $\widehat{\pi}^{(x,t)}$ is the piecewise linear function $\widehat{\pi}^{(x,t)} : (-\infty, t] \to \mathbb R$
with  $\widehat{\pi}^{(x,t)}(t-k) = \widehat{h}^{k}(x,t)(1)$ for every integer $k \geq 0$.
Let $\widehat{{\mathcal X}} := \{\widehat{\pi}^{(x,t)}: (x,t)\in \Z^2_{\text{odd}}\}$
be the collection of all possible dual paths admitted by $\widehat{G}$.
For the backward path $ \widehat{\pi} $, the scaled version is 
$ \widehat{\pi}_n : [  -\infty, \sigma_{\widehat{\pi}}/n] \to
[-\infty, \infty]$ given by $\widehat{\pi}_n (t)=
\widehat{\pi }(n t)/ \sqrt{n}$ for each $n \geq 1$. 
Let $ \widehat{{\cal X}}_n := \{\widehat{\pi}_n^{(x,t)} :(x,t) \in \Z^2_{\text{odd}}\}$ be the collection  of all the $n$ th order diffusively scaled dual paths. 
The following theorem regarding the joint 
convergence of $(\bar{{\cal X}}_n, \overline{\widehat{{\cal X}}}_n)$ to $(\WW, \widehat{\WW})$ is due to Roy \textit{ et. al.} \cite{RSS16}.
\begin{theorem}[Theorem 6.1 of \cite{RSS16}]
\label{th:SSRW_BW}
As $n\rightarrow \infty$,  $(\bar{{\cal X}}_n, \overline{\widehat{{\cal X}}}_n)$
 converges in distribution to the double Brownian web $({\cal W}, \widehat{\WW})$ as $({\cal H}\times \widehat{{\cal H}},{\cal B}_{{\cal H}\times \widehat{{\cal H}}})$ valued random variables.
\end{theorem}

We will use Theorem \ref{th:SSRW_BW}
to prove Theorem \ref{thm:ScalingLimCRTJt}. First we introduce 
some notation. For $t \in \R$ and $K \subseteq \Pi$, let $K^{t-} := \{\pi \in K : \sigma_{\pi} \leq t\}$ denote the collection of paths in $K$ which start before time $t$.
Similarly for $\widehat{K} \subseteq \widehat{\Pi}$, let $\widehat{K}^{t+} := \{\widehat{\pi} \in \widehat{K} : \sigma_{\widehat{\pi}} \geq t\}$ denote the collection of dual paths in $\widehat{K}$ which start after time $t$. 
Let $\xi_K \subset \R$ be defined as
 \begin{align}
  \label{def:EtaSet}
  \xi_K := \{ \pi(0) : \pi \in K^{(-1)-}\text{ with }\pi(0) \in [0,1]\}.
 \end{align}
For $x \in \xi_{K}$, for $K \subset \Pi$ and $\widehat{K}\subset \widehat{\Pi}$,
 let $\widehat{\pi}^{(x,0)}_r = \widehat{\pi}^{(x,0)}_r(K \times \widehat{K})$ be defined to be the path  $\widehat{\pi} \in \widehat{K}$ with $\sigma_{\widehat{\pi}} = 0 $ and such that there is no other path 
$\widehat{\pi}^\prime \in \widehat{K}^{0+} $ with $x< \widehat{\pi}^\prime(0)<\widehat{\pi}(0)$; if no such  
$\widehat{\pi} \in \widehat{K}$ exists then we define $\widehat{\pi}^{(x,0)}_r $ to be the {\it backward} constant zero function $\widehat{o}$ starting  at $\sigma_{\widehat{o}} = 0$.
In other words $\widehat{\pi}^{(x,0)}_r$ is the  path of $\widehat{K}^{0+}$ which intersects the $x$-axis to the right of $(x,0)$ and is the closest such path to do so.
Similarly we define $\widehat{\pi}^{(x,0)}_l$  as the path of $\widehat{K}^{0+}$ which intersects the $x$-axis to the left of $(x,0)$ and is the closest such path to do so. Let $\widehat{\gamma}(l,r) := \widehat{\gamma}(\widehat{\pi}^{(x,0)}_r, \widehat{\pi}^{(x,0)}_l)$  be the first meeting time  of the two backward paths $\widehat{\pi}^{(x,0)}_l$ and 
$\widehat{\pi}^{(x,0)}_r$ (which  can possibly be $- \infty$).
Let $\Delta = \Delta_{r,l}(x) \subset \R^2$ denote the region enclosed between these two backward paths, $\widehat{\pi}^{(x,0)}_l$ and $\widehat{\pi}^{(x,0)}_r$, and formally defined as 
\begin{align}
\label{def:Delta_x}
\Delta_{r,l}(x) := \{(y,s) : s \in (\widehat{\gamma}(l,r), 0), y \in (\widehat{\pi}^{(x,0)}_l(s), \widehat{\pi}^{(x,0)}_r(s)) \}.
\end{align}
 If $\widehat{\gamma}(l,r)$ is finite then this region is also bounded .
 
 Let ${\cal D} := \{(x/\sqrt{n},t/n) : (x,t) \in \Z^2_{\text{even}}, n \in \N\}$
 and 
\begin{align*}
{\cal S}(\widehat{\pi}^{(x,0)}_r, \widehat{\pi}^{(x,0)}_l) = 
{\cal S}(\widehat{\pi}^{(x,0)}_r, \widehat{\pi}^{(x,0)}_l)_{(K,\widehat{K})} & := \{\pi \in K : (\pi(\sigma_\pi),\sigma_{\pi}) \in \Delta_{r,l}(x) \cap {\cal D}\}, \\
\widehat{{\cal S}}(\widehat{\pi}^{(x,0)}_r, \widehat{\pi}^{(x,0)}_l) =
\widehat{{\cal S}}(\widehat{\pi}^{(x,0)}_r, \widehat{\pi}^{(x,0)}_l)_{(K,\widehat{K})} & := \{\widehat{\pi} \in \widehat{K} : 
(\widehat{\pi}(\sigma_{\widehat{\pi}}), \sigma_{\widehat{\pi}}) \in 
\Delta_{r,l}(x) \cap {\cal D}\},\\
\overline{{\cal S}}(\widehat{\pi}^{(x,0)}_r, \widehat{\pi}^{(x,0)}_l) = \overline{{\cal S}}(\widehat{\pi}^{(x,0)}_r, \widehat{\pi}^{(x,0)}_l)_{(K,\widehat{K})} & := \overline{{\cal S}(\widehat{\pi}^{(x,0)}_r, \widehat{\pi}^{(x,0)}_l)}_{(K,\widehat{K})} \text{ and } \\
\overline{\widehat{{\cal S}}}(\widehat{\pi}^{(x,0)}_r, \widehat{\pi}^{(x,0)}_l) = \overline{\widehat{{\cal S}}}(\widehat{\pi}^{(x,0)}_r, \widehat{\pi}^{(x,0)}_l)_{(K,\widehat{K})} & := \overline{\widehat{{\cal S}}(\widehat{\pi}^{(x,0)}_r, \widehat{\pi}^{(x,0)}_l)}_{(K,\widehat{K})}{\color{green} :}
\end{align*}
 the closures being taken in $(\Pi, d_\Pi)$ and $(\widehat{\Pi}, d_{\widehat{\Pi}})$ respectively. 
 We further assume that $(K,\widehat{K})$ is such that,
for each $x \in \xi_K$ both the collections, ${\cal S}(\widehat{\pi}^{(x,0)}_r, \widehat{\pi}^{(x,0)}_l)$ and $\widehat{{\cal S}}(\widehat{\pi}^{(x,0)}_r, \widehat{\pi}^{(x,0)}_l)$, are tree-like in the sense of Definition \ref{def:TreeLikeMetricFrmCoalPaths} and by construction both ${\cal S}(\widehat{\pi}^{(x,0)}_r, \widehat{\pi}^{(x,0)}_l)$ and $\widehat{{\cal S}}(\widehat{\pi}^{(x,0)}_r, \widehat{\pi}^{(x,0)}_l)$ are families of coalescing paths. 
 The completions of these spaces with respect to the ancestor metric are defined as 
\begin{align*}
{\cal M}^{(x,0)} = {\cal M}^{(x,0)}_{(K,\widehat{K})} & := {\cal M}({\cal S}(\widehat{\pi}^{(x,0)}_r, \widehat{\pi}^{(x,0)}_l)) \text{ and } \\
\widehat{{\cal M}}^{(x,0)} = \widehat{{\cal M}}^{(x,0)}_{(K,\widehat{K})} & := {\cal M}(\widehat{{\cal S}}(\widehat{\pi}^{(x,0)}_r, \widehat{\pi}^{(x,0)}_l)). 
\end{align*}
We further assume that for each $x \in \xi_K$, both these metric spaces, ${\cal M}^{(x,0)}$ and $\widehat{\cal M}^{(x,0)}$, are compact
and  let $\phi({\cal M}^{(x,0)})$ and $\phi^\prime(\widehat{\cal M}^{(x,0)})$ 
denote the isometric embeddings of the metric spaces ${\cal M}^{(x,0)}$ and $\widehat{\cal M}^{(x,0)}$ into $\mathbb{M}$ respectively, the metric space of isometric equivalence classes of compact
metric spaces endowed with the Gromov-Hausdorff metric.
Fix $f : \mathbb{M} \times \mathbb{M} \to \R$, a bounded continuous real valued function and define 
\begin{align*}
\kappa_{ (K, \widehat{K}) }(f):= 
\begin{cases}
0 & \text{ if } \xi_{K} = \emptyset\\
\sum_{x \in \xi_{K}\cap[0,1]} f(\phi({\cal M}^{(x,0)}), \phi^\prime(\widehat{\cal M}^{(x,0)})) &\text{ if }0 < \# \xi_{K} < \infty \\
\infty & \text{ otherwise}.
\end{cases}
\end{align*}

For the rest of the section we make two specific choices of $(K, \widehat{K})$, namely, $(\WW, \widehat{\WW})$ and $({\cal X}_n, \widehat{{\cal X}}_n)$. Also to simplify notation we write
$$
 \xi:= \xi_{\WW}, \quad  \xi_n :=\xi_{{\cal X}_n}, \quad \eta := \# \xi \text{ and }\eta_n := \# \xi_n' $$
It is known that $\E(\eta) = 1/\sqrt{\pi}$.

 As observed in \cite{RSS15},
 for each $x \in \xi$, there exist two dual paths
$\widehat{\pi}^{(x, 0)}_r, \widehat{\pi}^{(x, 0)}_l \in \widehat{\WW}$
both starting at $ (x, 0)$ such that $\widehat{\pi}^{(x, 0)}_r(-1) >  \widehat{\pi}^{(x, 0)}_l(-1) $ with $\widehat{\gamma}(r,l)$ finite. Hence the region $ \Delta_{r,l}(x)_{(\WW, \widehat{\WW})} $ is almost surely non-empty and bounded.
 From the properties of $(\WW, \widehat{\WW})$ we have that both the path spaces, ${\cal S}(\widehat{\pi}^{(x, 0)}_r,\widehat{\pi}^{(x, 0)}_l)_{(\WW \times \widehat{\WW})}$ and 
$\widehat{\cal S}(\widehat{\pi}^{(x, 0)}_r,\widehat{\pi}^{(x, 0)}_l)_{(\WW \times \widehat{\WW})}$, are tree-like in the sense of Definition
\ref{def:TreeLikeMetricFrmCoalPaths}.
Same argument as in Proposition \ref{prop:CRTCompact} shows 
that both the metric spaces, ${\cal M}^{(x,0)}_{(\WW, \widehat{\WW})}$
and $\widehat{{\cal M}}^{(x,0)}_{(\WW, \widehat{\WW})}$, are compact as well. Similarly 
for  $n \geq 1$ and for each $x_n \in \xi_n$, the region 
$\Delta_{r,l}(x)_{({\cal X}_n, \widehat{{\cal X}}_n)}$ is almost surely nonempty and bounded.
From the construction it follows that, both the families, ${\cal S}(\widehat{\pi}^{(x_n,0)}_r, \widehat{\pi}^{(x_n,0)}_l)_{({\cal X}_n , \widehat{{\cal X}}_n)}$ and $\widehat{\cal S}(\widehat{\pi}^{(x_n,0)}_r, \widehat{\pi}^{(x_n,0)}_l)_{({\cal X}_n, \widehat{{\cal X}}_n)}$, satisfy the condition
 of Definition \ref{def:TreeLikeMetricFrmCoalPaths} and both these metric spaces, ${\cal M}^{(x_n, 0)}_{({\cal X}_n, \widehat{{\cal X}}_n)}$ and $\widehat{{\cal M}}^{(x_n, 0)}_{({\cal X}_n, \widehat{{\cal X}}_n)}$, are compact.

Hence the random variables $\kappa(f) := \kappa_{({\mathcal  W},\widehat{\mathcal  W})}(f),
\text{ and } \kappa_n(f) := \kappa_{({\mathcal  X}_n,\widehat{{\mathcal  X}}_n)}(f) $ are well defined almost surely. 
In the next subsection we calculate $\E(\kappa(f))$.

\subsection{Calculation of $\E(\kappa(f))$}
\label{subsec:CalcKappa}

The aim of this subsection is to prove the following proposition which calculates $\E(\kappa(f))$.
\begin{proposition}
\label{prop:CalcKappaf}
We have
$ \E(\kappa(f)) = \E ( f(\phi(\TT), \phi^\prime(\widehat{\TT})))/\sqrt{\pi}$.
\end{proposition}

To prove this proposition we need to prove a limit result
for conditional continuum random trees. Recall that
 $\widehat{W}^{(0,0)}$ and $\widehat{W}^{(1/n,0)}$, are two independent coalescing backward Brownian motions starting from the 
points $(0,0)$ and $(1/n,0)$ respectively with their coalescing time $\widehat{\gamma}_n = \widehat{\gamma}(\widehat{W}^{(1/n,0)},\widehat{W}^{(0,0)})$.
 Let $\Delta(\widehat{W}^{(0,0)}, \widehat{W}^{(1/n,0)})$
denote the random region enclosed between these two backward Brownian paths, defined as
$$
\Delta( \widehat{W}^{(1/n,0)}, \widehat{W}^{(0,0)}) := 
\{(x,t) : t \in (\widehat{\gamma}_n, 0), x \in (\widehat{W}^{(0,0)}(t), \widehat{W}^{(1/n,0)}(t))\}.
$$

We construct an equivalent of $(\TT, \widehat{\TT})$ here.
 We start coalescing backward Brownian motions 
from all  points of $\Delta( \widehat{W}^{(1/n,0)}, \widehat{W}^{(0,0)}) \cap \Q^2$
such that these backwards paths coalesce with the boundary of 
$\Delta( \widehat{W}^{(1/n,0)}, \widehat{W}^{(0,0)})$ as soon as they intersect
 the boundary. This collection of coalescing backward paths is denoted by $\widehat{\cal S}( \widehat{W}^{(1/n,0)}, \widehat{W}^{(0,0)})$. 
From  earlier discussions it follows that for each $(x,t) \in \Delta( \widehat{W}^{(1/n,0)}, \widehat{W}^{(0,0)})\cap \Q^2$, there exists a forward path $\pi^{(x,t)} \in \Pi$ which does not cross any dual   path in $\widehat{\cal S}( \widehat{W}^{(1/n,0)}, \widehat{W}^{(0,0)})$ and defined as  
\begin{align}
\label{def:FwdCoalescBpathsinWn}
\pi^{(x,t)}(s) :=  
\begin{cases}
 \sup\{ y: y \in  \Q, \widehat{\pi}^{(y,s)} \in \widehat{{\cal S}}(\widehat{W}^{(1/n,0)}, \widehat{W}^{(0,0)}), \widehat{\pi}^{(y,s)}(t) < x\} = &\\
\quad = \inf\{y: y \in \Q,\widehat{\pi}^{(y,s)} \in \widehat{{\cal S}}(\widehat{W}^{(1/n,0)},\widehat{W}^{(0,0)}), \widehat{\pi}^{(y,s)}(t) > x\}&\text{ for }s \in [t,0) \cap \Q\\
\pi^{(x,t)}(0) &\text{ for }s\geq 0.
\end{cases}
\end{align}
We denote this collection of forward paths starting from all the points of  $\Delta( \widehat{W}^{(1/n,0)}, \widehat{W}^{(0,0)})\cap \Q^2$ by ${\cal S}( \widehat{W}^{(1/n,0)}, \widehat{W}^{(0,0)})$. Because of property (d) of $(\WW, \widehat{\WW})$, the set $\{\pi(0) : \pi \in {\cal S}(\widehat{W}^{(1/n,0)}, \widehat{W}^{(0,0)})$, $\sigma_{\pi} \leq -1\}$ is almost surely finite.
Let $x^n = x^n(\omega) \in (0,1/n)$ be defined as 
\begin{align*}
x^n := \min\{x \in (0,1/n) : \text{ there exists }\pi \in {\cal S}(\widehat{W}^{(1/n,0)}, \widehat{W}^{(0,0)})\text{ with }\sigma_\pi \leq -1\text{ and } \pi(0) = x\}.
\end{align*}
Next we consider all the paths in ${\cal S}( \widehat{W}^{(1/n,0)}, \widehat{W}^{(0,0)})$
which pass through $x^n$. In other words we consider a new collection of coalescing paths ${\cal S}^\prime(\widehat{W}^{(1/n,0)}, \widehat{W}^{(0,0)})$ defined as ${\cal S}^\prime(\widehat{W}^{(1/n,0)}, \widehat{W}^{(0,0)}):= \{\pi \in {\cal S}( \widehat{W}^{(1/n,0)}, \widehat{W}^{(0,0)}) : \pi(0) = x^n\}$.
By construction, both ${\cal S}^\prime( \widehat{W}^{(1/n,0)}, \widehat{W}^{(0,0)})$ and $\widehat{\cal S}( \widehat{W}^{(1/n,0)}, \widehat{W}^{(0,0)})$
are tree-like in the sense of Definition \ref{def:TreeLikeMetricFrmCoalPaths}. 
For $n \geq 1$ we consider the pair of complete metric spaces
 $$
 (\TT_n, \widehat{\TT}_n) := \bigl({\cal M}({\cal S}^\prime( \widehat{W}^{(1/n,0)}, \widehat{W}^{(0,0)})),{\cal M}(\widehat{\cal S}( \widehat{W}^{(1/n,0)}, \widehat{W}^{(0,0)}))\bigr).
 $$ 
 Similar argument as in the proof of Proposition \ref{prop:CRTCompact}
shows that these metric spaces are also compact almost surely. It is useful to observe that for each $n \geq 1$, the metric space $\TT_n$ does not depend on how the forward paths evolve after time 
$0$.
 
 \begin{lemma}
 \label{lem:CRTConditionalLimit}
  $$(\TT_n, \widehat{\TT}_n) | \{\widehat{\gamma}_n < -1\} \Rightarrow (\TT, \widehat{\TT})\text{ as }n\to \infty,$$
  where convergence in distribution holds with respect to the Gromov-Hausdorff topology.
 \end{lemma}
\noindent Assuming the above lemma we first calculate $\E(\kappa(f))$ and complete the proof 
of Proposition \ref{prop:CalcKappaf}.

\noindent {\bf Proof of Proposition \ref{prop:CalcKappaf} :}
Let $I_n \subset \{0,1,\ldots,n-1\}$ be given by $$I_n :=
\{i : 0 \leq i \leq n-1\text{ such that } \widehat{\pi}^{(i/n,0)} , \widehat{\pi}^{((i+1)/n,0)} \in \widehat{\mathcal  W} \text{ with }\widehat{\pi}^{(i/n,0)}(-1) <
\widehat{\pi}^{((i+1)/n,0)}(-1)\}.$$ For $i \in I_n$, let 
$\Delta(\widehat{\pi}^{((i+1)/n,0)}, \widehat{\pi}^{(i/n,0)})$
denote the random region enclosed by the two backward Brownian paths $\widehat{\pi}^{((i+1)/n,0)}$ and $\widehat{\pi}^{(i/n,0)}$ and let
\begin{align*}
{\cal S}(\widehat{\pi}^{((i+1)/n,0)}, \widehat{\pi}^{(i/n,0)})
& := \{\pi \in \WW : (\pi(\sigma_\pi),\sigma_\pi)\in \Delta(\widehat{\pi}^{((i+1)/n,0)}, \widehat{\pi}^{(i/n,0)})\cap \Q^2\} \text{ and }\\
\widehat{\cal S}(\widehat{\pi}^{((i+1)/n,0)}, \widehat{\pi}^{(i/n,0)})
& := \{\widehat{\pi} \in \widehat{\WW} : (\widehat{\pi}(\sigma_{\widehat{\pi}}),\sigma_{\widehat{\pi}})\in \Delta(\widehat{\pi}^{((i+1)/n,0)}, \widehat{\pi}^{(i/n,0)})\cap \Q^2\}.
\end{align*} 
For $i \in I_n$, let $x^n_i = x^n_i(\omega) \in (i/n,(i+1)/n)$ be defined as 
\begin{align*}
x^n_i := \min\{x \in (i/n,(i+1)/n) : \text{ there exists }\pi \in {\cal S}(\widehat{\pi}^{((i+1)/n,0)}, \widehat{\pi}^{(i/n,0)})\text{ with }\sigma_\pi\leq - 1\text{ and } \pi(0) = x\}.
\end{align*}
We recall that for each $i \in I_n$, the set $\{\pi(0) : \pi \in \WW^{(-1)-}, \pi(0) \in (i/n , (i+1)/n)\}$ is nonempty and finite and hence $x^n_i$ is well defined.
Next we define the coalescing family of forward paths passing through $x^n_i$ given as 
$${\cal S}^\prime(\widehat{\pi}^{((i+1)/n,0)}, \widehat{\pi}^{(i/n,0)}) := 
\{\pi \in {\cal S}(\widehat{\pi}^{((i+1)/n,0)}, \widehat{\pi}^{(i/n,0)}): \pi (0) = x^n_i \}.$$
 From the properties of $(\WW, \widehat{\WW})$ it follows that both these collections of 
paths, ${\cal S}^\prime(\widehat{\pi}^{((i+1)/n,0)}, \widehat{\pi}^{(i/n,0)})$ and 
$\widehat{\cal S}(\widehat{\pi}^{((i+1)/n,0)}, \widehat{\pi}^{(i/n,0)})$, are tree-like
and we consider the complete metric spaces 
$$
{\cal M}({\cal S}^\prime(\widehat{\pi}^{((i+1)/n,0)}, \widehat{\pi}^{(i/n,0)}))
\text{ and }
{\cal M}(\widehat{\cal S}(\widehat{\pi}^{((i+1)/n,0)}, \widehat{\pi}^{(i/n,0)})).
$$
For $i \in I_n$, the paths in ${\cal S}^\prime(\widehat{\pi}^{((i+1)/n,0)}, \widehat{\pi}^{(i/n,0)})$
do not cross the dual paths in $\widehat{\cal S}(\widehat{\pi}^{((i+1)/n,0)}, \widehat{\pi}^{(i/n,0)})$ and $(\widehat{\pi}^{((i+1)/n,0)}, \widehat{\pi}^{(i/n,0)})$ is
distributed as coalescing backward Brownian motions starting from $((i+1)/n,0)$ and $(i/n, 0)$ conditioned to meet after time $-1$. Hence for each $i \in I_n$ we have
\begin{align}
\label{eq:KappaDistEq}
(\TT_n,\widehat{\TT}_n)|\{\widehat{\gamma}_n < -1\}
\stackrel{d}{=} \bigl({\cal M}({\cal S}^\prime(\widehat{\pi}^{((i+1)/n,0)}, \widehat{\pi}^{(i/n,0)})),{\cal M}(\widehat{\cal S}(\widehat{\pi}^{((i+1)/n,0)}, \widehat{\pi}^{(i/n,0)}))\bigr).
\end{align}
Fix a bounded continuous function $g : \mathbb{M}\times \mathbb{M} \to \R$.
We define
$${\cal R}_n(g) := \sum_{i \in I_n} g\Bigl(\phi_{i/n}\bigl({\cal M}({\cal S}^\prime(\widehat{\pi}^{((i+1)/n,0)}, \widehat{\pi}^{(i/n,0)}))\bigr),\phi^\prime_{i/n}\bigl({\cal M}(\widehat{\cal S}(\widehat{\pi}^{((i+1)/n,0)}, \widehat{\pi}^{(i/n,0)}))\bigr)\Bigr),$$
where $\phi_{i/n}$ and $\phi^\prime_{i/n}$ are isometric embeddings of the respective metric spaces into $\mathbb{M}$, the space of all isometry equivalence classes of compact metric spaces.  
We want to show that $${\cal R}_n(g) \to \sum_{x \in \xi} g\bigl(\phi_x({\cal M}^{(x,0)}_{(\WW,\widehat{\WW})}),\phi^\prime_x(\widehat{\cal M}^{(x,0)}_{(\WW,\widehat{\WW})})\bigr)\text{ as }n \to \infty\text{ almost surely},$$
where $\phi_x, \phi^\prime_x$ are the corresponding embeddings of the respective metric spaces into $\mathbb{M}$.

We first show that $\# I_n \to \# \xi = \eta$ as $n \to \infty$. Choose $n_0 = n_0(\omega)$
such that for all $x,y\in \xi$ and for $n \geq n_0$ we have $|x-y| > 2/n$.
For each $x \in \xi$, set   $ l^x_n =\lfloor nx \rfloor/n $
and $ r^x_n  =  l^x_n + (1/n)$.
From the non-crossing nature of paths in $\widehat{\WW}$ it follows that for each $x\in \xi$, we have
$\widehat{\pi}^{(r^x_n,0)}(-1) > \widehat{\pi}^{(l^x_n,0)}(-1)$ and hence $l^x_n \in I_n$.
 Hence $\# I_n \geq \xi$ for all $n\geq n_0$.
On the other hand for each $i \in I_n$, choose $y_i$ such that 
$y_i \in (\widehat{\pi}^{(i/n,0)}(-1),
\widehat{\pi}^{((i+1)/n,0)}(-1))\cap \Q$ and consider the path $\pi^{(y_i,-1)}$ in ${\mathcal  W}$. Then we must have
$\pi^{(y_i,-1)}(0) = x_i \in \xi$. This implies that 
 $\# I_n \leq \# \xi$ for all $n$. 
Hence  for all $n\geq n_0$ we have $\# I_n =\# \xi$.

This ensures that for all large $n$,  
\begin{equation}
\label{eq:RntogM1}
{\cal R}_n(g) = \sum_{x \in \xi} g\Bigl(\phi_{l^x_n}\bigl({\cal M}({\cal S}^\prime(\widehat{\pi}^{(r^x_n,0)}, \widehat{\pi}^{(l^x_n,0)}))\bigr),\phi^\prime_{l^x_n}\bigl({\cal M}(\widehat{\cal S}(\widehat{\pi}^{(r^x_n,0)}, \widehat{\pi}^{(l^x_n,0)}))\bigr)\Bigr).
\end{equation}

Choose $n_1 = n_1(\omega)\geq n_0$ large enough such that 
for all $n \geq n_1$ we have $\# I_n = \eta$ and consequently for all $n \ge n_1$ and for each $x \in \xi$ we have
\begin{align}
\label{eq:RntogM2}
{\cal M}({\cal S}^\prime(\widehat{\pi}^{(r^x_n,0)}, \widehat{\pi}^{(l^x_n,0)}))
= {\cal M}^{(x,0)}_{(\WW,\widehat{\WW})}.
\end{align}

We recall that there are exactly two dual paths $\widehat{\pi}^{(x,0)}_r$ and $\widehat{\pi}^{(x,0)}_l$ starting from $(x,0)$ with $\widehat{\pi}^{(x,0)}_r(-1) > \widehat{\pi}^{(x,0)}_l(-1)$. Hence
from property (f) of $({\cal W},\widehat{\cal W})$ it follows that, 
$
\{\widehat{\pi}^{(l^x_n,0)} : n \in \N\}$ and $\{\widehat{\pi}^{(r^x_n,0)} : n \in \N\}$
converge to $\widehat{\pi}^{(x,0)}_l$ and $\widehat{\pi}^{(x,0)}_r$ respectively
in $(\widehat{\Pi}, d_{\widehat{\Pi}})$ as $n \to \infty$.
Let $\widehat{\gamma}(\widehat{\pi}^{(x,0)}_l,\widehat{\pi}^{(l^x_n,0)})$ denote the coalescing time of the two dual paths
$\widehat{\pi}^{(l^x_n,0)}$ and $\widehat{\pi}^{(x,0)}_l$. 
Similarly the coalescing time of the two dual paths
$\widehat{\pi}^{(r^x_n,0)}$ and $\widehat{\pi}^{(x,0)}_r$ is given by 
$\widehat{\gamma}(\widehat{\pi}^{(x,0)}_r,\widehat{\pi}^{(r^x_n,0)})$.
For all $n \ge n_0$ and for each $x \in \xi$, the metric space  $\widehat{{\cal M}}^{(x,0)}_{(\WW,\widehat{\WW})}$ is naturally embedded into ${\cal M}(\widehat{\cal S}(\widehat{\pi}^{(r^x_n,0)}, \widehat{\pi}^{(l^x_n,0)}))$. Using this embedding we obtain that for all $n \ge n_1$ and for each $x \in \xi$ we have
\begin{align}
\label{eq:RntogM3}
d_{\text{GH}}\bigl({\cal M}(\widehat{\cal S}^\prime(\widehat{\pi}^{(r^x_n,0)}, \widehat{\pi}^{(l^x_n,0)})),\widehat{\cal M}^{(x,0)}_{(\WW,\widehat{\WW})} \bigr)\leq 
2(\widehat{\gamma}(\widehat{\pi}^{(x,0)}_l,\widehat{\pi}^{(l^x_n,0)})\vee
\widehat{\gamma}(\widehat{\pi}^{(x,0)}_r,\widehat{\pi}^{(r^x_n,0)})).
\end{align}
Now because of property (g) of $(\WW, \widehat{\WW})$, we have that $(\widehat{\gamma}(\widehat{\pi}^{(x,0)}_l,\widehat{\pi}^{(l^x_n,0)}),\widehat{\gamma}(\widehat{\pi}^{(x,0)}_r,\widehat{\pi}^{(r^x_n,0)}))  \to (0,0)$ as $n \to \infty$. 
Since $g$ is continuous on the product metric space and the set $\xi$ is finite, (\ref{eq:RntogM1}), (\ref{eq:RntogM2}) and (\ref{eq:RntogM3}) show that 
$${\cal R}_n(g) \to \sum_{x \in \xi} g\bigl(\phi_x({\cal M}^{(x,0)}_{(\WW,\widehat{\WW})}),\phi^\prime_x(\widehat{\cal M}^{(x,0)}_{(\WW,\widehat{\WW})})\bigr)\text{ as }n \to \infty.$$
As $g$ is bounded and $\E[\eta] < \infty$, the family $\{{\cal R}_n(g) : n \in \N \}$ is uniformly integrable and hence we have 
$$
\lim_{n \to \infty} \E[{\cal R}_n(g)] = \E\bigl(\sum_{x \in \xi} g\bigl(\phi_x({\cal M}^{(x,0)}_{(\WW,\widehat{\WW})}),\phi^\prime_x(\widehat{\cal M}^{(x,0)}_{(\WW,\widehat{\WW})})\bigr)\bigr).$$
From (\ref{eq:KappaDistEq}) and  for $B\in[0,\infty)$ denoting a standard Brownian motion we have,
\begin{align*}
& \lim_{n \to \infty} \E\bigl({\cal R}_n(g)\bigr)\\
& = \lim_{n \to \infty}n \E\bigl(g(\TT_n,\widehat{\TT}_n) \bigl | 
 \widehat{\gamma}_n < -1\bigr)\P(\widehat{\gamma}_n < -1)\\
 & = \lim_{n \to \infty}n \E\bigl(g(\TT_n,\widehat{\TT}_n) \bigl | 
 \widehat{\gamma}_n < -1\bigr)\P(1/n + \min_{t\in [0,1]}\sqrt{2} B(t) > 0) \\
& = \lim_{n \to \infty} \E\bigl(g(\TT_n,\widehat{\TT}_n) \bigl | 
 \widehat{\gamma}_n < -1\bigr) n (2\Phi(1/\sqrt{2}n) - 1) \nonumber \\
& =  \E \bigl( g(\TT,\widehat{\TT}) \bigr) /\sqrt{\pi}
\end{align*}
 where $\Phi$ is the distribution function of a standard normal random variable.
The last step follows from Lemma \ref{lem:CRTConditionalLimit}.
This completes the proof.
\qed

\noindent \textbf{ Proof of Lemma \ref{lem:CRTConditionalLimit} :}
Recall that $\widehat{B}_1$ and $ \widehat{B}_2$ are two independent
backward Brownian motions both starting from the origin.  We observe that $(\widehat{B}_1, \widehat{B}_2)\in A_n$  almost surely for all $n\geq 1$. 
The proof of Corollary \ref{cor:DeterministicContinuousMapping} shows that 
$(\lambda^n_1, \lambda^n_0) \to (\lambda_1, \lambda_0)$ as $n \to \infty$ almost surely.
Choose $n_0 = n_0(\omega)$ so that for all $n \geq n_0$ we have 
$\lambda^n_0 \in (\lambda_1, \lambda_0)$ and consequently $\lambda^n_1 = \lambda_1$.

For all $n \geq n_0$, the set $\{\pi(\lambda^n_0 - \lambda_0) : \pi \in {\cal S}(\widehat{B}^+, \widehat{B}^-), \sigma_{\pi} \leq (\lambda^n_0 - \lambda_0) - 1\}$ is nonempty and finite.
We take $y^n \in \R$ as $y^n := \min\{\pi(\lambda^n_0 - \lambda_0) : \pi \in {\cal S}(\widehat{B}^+, \widehat{B}^-), \sigma_{\pi} \leq (\lambda^n_0 - \lambda_0) - 1\}$.
Next we consider the coalescing families of forward and backward paths defined by
\begin{align*}
{\cal S}_n(\widehat{B}^+, \widehat{B}^-) & := \{\pi \in ({\cal S}_n(\widehat{B}^+, \widehat{B}^-))^{(\lambda^n_0 - \lambda_0)-}: \pi(\lambda^n_0 - \lambda_0) = y^n \}\text{ and }\\
\widehat{\cal S}_n(\widehat{B}^+, \widehat{B}^-) & := \{\widehat{\pi} \in \widehat{\cal S}_n(\widehat{B}^+, \widehat{B}^-): \sigma_{\widehat{\pi}} \leq (\lambda^n_0 - \lambda_0)\}.
\end{align*}
Let ${\cal M}({\cal S}_n(\widehat{B}^+, \widehat{B}^-))$ and ${\cal M}(\widehat{\cal S}_n(\widehat{B}^+, \widehat{B}^-))$ denote the corresponding complete metric spaces 
where completions are taken with respect to the ancestor metric.
We observe that both these metric spaces are naturally embedded into the metric spaces 
${\cal M}({\cal S}(\widehat{B}^+, \widehat{B}^-)) = \TT$ and ${\cal M}(\widehat{\cal S}(\widehat{B}^+, \widehat{B}^-)) = \widehat{\TT}$ respectively.
We show that almost surely
\begin{equation}
\label{eq:ConvMSntoMS}
d_{GH}({\cal M}({\cal S}_n(\widehat{B}^+, \widehat{B}^-)),\TT)\vee
d_{GH}({\cal M}(\widehat{\cal S}_n(\widehat{B}^+, \widehat{B}^-)),\widehat{\TT}) \to 0,
\end{equation}
as $n \to \infty$. Fix $\epsilon > 0$ and choose $n_1 = n_1(\omega)\geq n_0$ such that $( \lambda_0 - \lambda^n_0) < \epsilon/2$ for all $n \geq n_1$.
We first observe that the natural embedding shows that 
$d_{GH}({\cal M}(\widehat{\cal S}_n(\widehat{B}^+, \widehat{B}^-)),\widehat{\TT}) \leq 
2 ( \lambda_0 - \lambda^n_0)< \epsilon$ for all $n \geq n_1$. 

In order to get an upper bound for $d_{GH}({\cal M}({\cal S}_n(\widehat{B}^+, \widehat{B}^-)),{\TT})$, from Remark \ref{rem:CRTRootType} we obtain that for all $\pi_1, \pi_2 \in {\cal S}(\widehat{B}^+, \widehat{B}^-)$ with $\sigma_{\pi_1}\vee \sigma_{\pi_2} \leq -\epsilon/2$
there exists $s_0 = s_0(\omega) \in (0,\epsilon/2) $ which does not depend on the choice of $\pi_1,\pi_2$ such that  for all $s \geq -s_0$ we have $\pi_1(s) = \pi_2(s)$. Choose $n_2 = n_2(\omega)\geq n_1$
 such that $1/n_2 < s_0$.
Then for any $\pi \in {\cal S}(\widehat{B}^+, \widehat{B}^-)$ with $\sigma_{\pi}
 \leq -\epsilon/2$ we have $\pi \in {\cal S}_n(\widehat{B}^+, \widehat{B}^-)$ as well for all $n \geq n_2$. Hence we have $d_{GH}({\cal M}({\cal S}_n(\widehat{B}^+, \widehat{B}^-)),{\TT})
 \leq \epsilon$ for all $n \geq n_2$. This proves (\ref{eq:ConvMSntoMS}).

Finally same argument as in Lemma \ref{lem:ConditionalStoppingTime} shows that 
 $(\TT_n,\widehat{\TT}_n)|\{\widehat{\gamma}_n < - 1\}$ converges in distibution to $(\TT, \widehat{\TT})$ as $n\to \infty$.
 \qed

\subsection{ Convergence of $\E(\kappa_n(f))$}

In this subsection we prove that 
$\E(\kappa_n(f)) \to \E(\kappa(f))$ as $n \to \infty$
and use it to prove Theorem \ref{thm:ScalingLimCRTJt}. In order to prove the above mentioned convergence, we first need to show that $\E(\eta_n) \to \E(\eta)$ as $n \to \infty$.
For $t_1 < t_2$ let
\begin{align*}
\xi(t_1, t_2) & := \{\pi(t_2) : \pi \in \WW^{t_1-}, \pi(t_2) \in [0,1]\} \text{ and }\eta(t_1, t_2) := \#\xi(t_1, t_2),\\
\xi_n(t_1, t_2) & := \{\pi(t_2) : \pi \in {\cal X}_n^{t_1-}, \pi(t_2) \in [0,1]\} \text{ and }\eta_n(t_1, t_2) := \#\xi_n(t_1, t_2).
\end{align*}
 Below we prove a continuity property of $\eta(s,0) : s < 0$, which will be used in proving convergence of $\E(\eta_n)$. 

\begin{corollary}
\label{cor:ContinuityOccupdPtSet}
We have $\lim_{s \downarrow -1}\eta(s, 0) = \eta$ a.s.
\end{corollary}
\noindent \textbf{ Proof :} 
For $-1 \leq s_1 \leq s_2 < 0$ we have $\xi(s_1, 0) \subseteq 
\xi(s_2, 0)$ and consequently $\eta(s_1, 0) \leq \eta(s_2, 0)$. Because of this monotonicity $\lim_{s \downarrow -1}\eta(s, 0)$
exists almost surely. Since $ \xi \subseteq \xi(s,0)$ for all $s \in (-1,0)$,
we have $ \xi \subseteq \bigcap_{ s \in (-1,0)} \xi(s,0)$, implying that
$\lim_{s \downarrow -1} \eta(s,0) \geq \eta$.
If $\lim_{s \downarrow -1} \eta(s,0) 
> \eta$, then there exists  $x \in  \xi(-1/2 , 0 ) 
\setminus \xi$ and  
a sequence of paths $\{\pi^n \in {\cal W} \}$ with $\pi^n(0) = x$ for all $n$
such that $  -1 < \sigma_{\pi^n} \leq -1/2$ and  
$ \sigma_{\pi^n} $ decreases to $ -1 $ as $n \to \infty$.  
By compactness of ${\cal W}$, 
it follows that there exists a convergent subsequence $\{\pi^{n_k}: k \in \N\}$
and $\pi \in {\cal W}$ such that
$\pi^{n_k} \to \pi$ in $(\Pi, d_{\Pi})$ as $k \to \infty $. Since convergence in $(\Pi,d_\Pi)$ implies convergence of starting times also, we must have 
$\sigma_{\pi} = -1$ and $\pi(0) = x$ which 
contradicts the choice of $x$ and the proof follows.
\qed 

\noindent 
Next we use Corollary \ref{cor:ContinuityOccupdPtSet} to prove convergence of $\E(\eta_n)$.

\begin{lemma}
\label{lem:OccupdPtSetExpConv}
 $\E[\eta_n] \to \E[\eta ]$ as $n \to \infty$.
\end{lemma}
\noindent \textbf{ Proof :}
We show that $\eta_n$ converges in distribution
to $\eta$ as $n \to \infty$ and the sequence 
$\{\eta_n : n \in \N\}$ is uniformly integrable as well.
Recall from Theorem \cite{FINR04} that 
$(\bar{\mathcal X}_n,\bar{\widehat{\mathcal X}}_n) \Rightarrow 
({\mathcal W}, \widehat{{\mathcal W}})$ as $n \to \infty$. 
Using Skorohod's representation theorem  
we assume that we are working on a probability space such that
\begin{equation*}
d_{{\mathcal H}\times \widehat{{\mathcal H}}}((\bar{\mathcal X}_n, 
\bar{\widehat{\mathcal X}}_n), ({\mathcal W}, \widehat{{\mathcal W}})) \to 0,
\end{equation*}
almost surely as $ n \to \infty$. Same argument as in Lemma 3.2 of \cite{RSS15}
shows that $\eta_n$ converges to $\eta$ almost surely as $n \to \infty$. For completeness we present the proof here also.
First we show that, for all $k \geq 0$, 
\begin{equation}
\label{eqn:weakconvliminf}
\liminf_{n\to \infty} \mathbf{1}_{\{\eta_n \geq k\}} \geq \mathbf{1}_{\{\eta \geq k\}}
\text{ almost surely}.
\end{equation}
Indeed, for $k = 0$, both $\mathbf{1}_{\{\eta_n \geq k\}}$ and $\mathbf{1}_{\{\eta \geq k\}}$ equal $1$. We need to show (\ref{eqn:weakconvliminf}) for $k \geq 1$. 

From Properties of $(\WW, \widehat{\WW})$ it follows that 
$\eta$ is finite and $\xi \subset (0,1)$ almost surely. Hence  
we can choose $0 <\epsilon = \epsilon(\omega)$ such that 
\begin{itemize}
 \item[(i)] for all $x,y \in \xi$ with $x \neq y$ we have $(x -\epsilon, x +\epsilon) \subset (0,1)$ and $|x-y|> 2 \epsilon$;
\item[(ii)] for all $x \in \xi$ there exists $\pi \in {\cal W}$ with 
$\sigma_\pi \leq -1 -\epsilon$ and $\pi(0) = x$.
\end{itemize}
Let  $m = m(\omega) > 0$ be such that 
\begin{align}
\label{eq:dualCompact}
\{(y,t): \widehat{\pi}^{(0,0)}(t) \leq y \leq \widehat{\pi}^{(1,0)}(t), t \in [-1,0],
 \widehat{\pi}^{(0,0)}, \widehat{\pi}^{(1,0)} \in \widehat{{\cal W}}\}
\subset [-m, m]^2.
\end{align}
Choose $n_0 = n_0(\omega)$ such that for all $n \geq n_0$, 
$d_{{\cal H}\times \widehat{\cal H}}(({\cal W},\widehat{{\cal W}}),(\bar{{\cal X}}_n,
\bar{\widehat{{\cal X}}}_n)) < g(\epsilon,m)$ where $g(\epsilon,m)$ 
is as in Remark \ref{rem:MetricConv}.
By the choice of $n_0$, it follows that for each $\pi \in {\cal W}^{(-1-\epsilon)-}$ 
with $\pi(0) \in [0,1]$ and for all $n\geq n_0$, there 
exists $\pi^n \in {\cal X}_n$ with $\sigma_{\pi^n} \leq -1$ and $d_{\Pi}(\pi,\pi^n) < g(\epsilon,m)$. Since $|\pi(0) - \pi^n(0)| < \epsilon$, 
from the choice of $\epsilon$ we must have $\pi^n (0)
\in (0,1)$. Thus $\eta_n(\omega) \geq k$ for all $n \geq n_0$.
This proves (\ref{eqn:weakconvliminf}).

Finally we need to show that 
$\P(\limsup_{n \to \infty} \{\eta_n > \eta \}) = 0$.
This is equivalent to showing that $\P(\Omega^{k}_0) = 0$ for all $k\geq 0$,
where 
$$
\Omega^{k}_0 := \{\omega: \eta_n(\omega) > \eta(\omega) = k \text{ for infinitely 
many }n\}.
$$
Consider $k=0$ first. Using Corollary \ref{cor:ContinuityOccupdPtSet}
we can obtain $l_0:=l_0(\omega) \in (-1,0)$
such that $\eta(l_0,0)(\omega) = \eta(\omega)$.
Since both the points $(0,0)$ and $(1,0)$ are of type $(0,1)$ and 
on the event $\{\eta = 0\}$ 
the set $\xi(l_0,0) = \xi = \emptyset$,
we can obtain $\epsilon:=\epsilon(\omega) \in (0, 1+l_0)$ 
such that for all $\pi \in {\cal W}(\omega)$
with $\sigma_{\pi} \leq -1 + \epsilon \leq l_0$, $\pi(0) \notin (-\epsilon,1+\epsilon)$. 
On the event $\Omega^0_0$, choose large enough $n$ such that $\eta_n(\omega) > 0$ and 
$d_{{\cal H}\times \widehat{\cal H}}(({\cal W},\widehat{{\cal W}}),(\bar{{\cal X}}_n,
\bar{\widehat{{\cal X}}}_n)) < g(\epsilon,m)$ where $m = m(\omega) > 0$ is chosen as in 
(\ref{eq:dualCompact}) and $g(\epsilon,m)$ is as in Remark \ref{rem:MetricConv}.
From the choice of $\epsilon$, it follows that there exists $\pi \in {\cal W}(\omega)$
with $\sigma_{\pi} \leq -1 + \epsilon$ and $\pi(0) \in (-\epsilon, 1 + \epsilon)$. This 
gives a contradiction. Hence we have $\P(\Omega^0_0) = 0$.

For $k > 0$, on the event $\Omega^{k}_0$ we show that
a forward path $\pi \in {\mathcal W}$ coincides with a dual path $\widehat{\pi}
\in \widehat{\mathcal W}$ for a positive time which leads to a contradiction. 
From Corollary \ref{cor:ContinuityOccupdPtSet} it follows that,
we can choose $-1<l_0(\omega)< - 1/2 $ such that $\xi(l_0,0)(\omega) = \xi(\omega)$.
From property (a) of $({\cal W},\widehat{\cal W})$ it further follows that
for any $x \in \xi$, the type of the point $(x,0)$ is $(1,1)$. 
Hence for each $x \in \xi$ there exists $s_x = s_x(\omega)\in (-1,0)$ such that for any $\pi^1,\pi^2 \in {\cal W}^{l_0 -}$ with
$\pi^1(0) = \pi^2(0) = x$, we have $\pi^1(s) = \pi^2(s)$ for all $s \geq s_x$.
Set $s_0 := \max\{s_x : x \in \xi \}$ and observe that $s_0  < 0 $ almost surely.
By definition, for all $s \geq s_0 $ we have $\# \{\pi(s_0) : \pi \in {\cal W}^{l_0-}, \pi(1) \in [0,1]\} = k$. This gives us that $\xi(l_0,0)(\omega) = \xi(\omega) = \xi(l_0,s_0)(\omega)$, i.e., the paths leading to any single point considered in $ \xi(l_0,0)$
 have coalesced before time $ s_0$. We consider $k$ paths contributing to $\xi(l_0, 0) $, {\it viz}\/., $\pi_1, \dotsc, \pi_k$ 
such that $\{\pi_1(0), \dotsc, \pi_k(0)\} = \xi(l_0, 0)$. Choose $0 <\epsilon = 
\epsilon(\omega) < \min\{|s_0|, 1 + l_0 \}/3$ such that
\begin{itemize}
 \item[(i)] $(x - \epsilon,x + \epsilon) \subset (0,1)$ for all $x \in \xi(l_0, 0)$;
 \item[(ii)]the $\epsilon$-tubes around $\pi_1 , \dotsc, \pi_k$, given by 
$$
T^{i}_{\epsilon} := \{(x,t): |x-\pi_i(t)| < \epsilon ,
s_0 \leq t \leq 0 \}\text{ for }i=1,\dotsc,k
$$
are disjoint.
\end{itemize}
Choose $m = m(\omega) > 0$ as in (\ref{eq:dualCompact}). 
Set $n_0 = n_0(\omega)$ such that for all $n \geq n_0$, $\xi_{n_0}(\omega) > k$
and  
$$d_{{\cal H}\times \widehat{\cal H}}(({\cal W},\widehat{{\cal W}})(\omega),(\bar{{\cal X}}_n,
\bar{\widehat{{\cal X}}}_n)(\omega)) < g(\epsilon,m)$$ where $g(\epsilon,m)$ 
is as in Remark \ref{rem:MetricConv}.
By the choice of $n_0$, it follows that one of the $k$ tubes must contain at least two paths, $ \pi_1^{n_0}, \pi_2^{n_0}$ (say) 
of ${\mathcal X}^{(-1)-}_{n_0}$ which do not coalesce by time $0$.
From the construction of dual paths it follows that there exist at least one dual
path $\widehat{\pi}^{n_0} \in \widehat{{\mathcal X}}_{n_0}$ 
lying between $\pi_{1}^{n_0}$ and $\pi_{2}^{n_0}$ for $t\in [-1,0]$ and hence 
we must have an approximating $\widehat{\pi} \in \widehat{\mathcal W}$ close to 
$\widehat{\pi}^{n_0}$ for $t \in [-1,-\epsilon]$. More formally there exists a dual path $\widehat{\pi} \in \widehat{\cal W}^{(-\epsilon)+}$ such that $\sup_{t \in [-1,  - \epsilon]}|\widehat{\pi}(t) - \widehat{\pi}^{n_0}(t)|< \epsilon$ and thus $\sup_{t \in [s_0, - \epsilon]}
|\widehat{\pi}(t)- \pi(t)|< 4\epsilon $.

Now select any sequence $ \epsilon_l \downarrow 0 $ with $ \epsilon_l < \epsilon $ for 
all $l \geq 1$. By the above procedure we can select $ \widehat{\pi}_l \in  
\widehat{\cal W}$ with $\sigma_{\widehat{\pi}_l} \geq - \epsilon_l > - \epsilon$
such that $ \sup_{t \in [s_0,-\epsilon]}|\widehat{\pi}_l(t) - \pi(t)|
< 4\epsilon_l $ for some $ \pi \in \{ \pi_1, \dotsc, \pi_k \}$. Since we have only finitely many paths $ \pi_1, \dotsc, \pi_k$,  we can choose $ \pi_i $ for some $ 1 \leq i \leq k$ and a subsequence $ \widehat{\pi}_{l_j} $ so that $  \sup_{t \in [s_0,-\epsilon]}|\widehat{\pi}_{l_j}
(t) - \pi_i(t)|< 4\epsilon_{l_j}  $ for all $ j \geq 1$.
By the compactness of $  \widehat{\cal W} $,
there exists $ \widehat{ \pi} \in \widehat{\cal W} $ with  $\sigma_{\widehat{\pi}}
\geq - \epsilon$ such that $ \widehat{\pi} (t) =  \pi_i(t)
$ for $ t \in [s_0, - \epsilon]$.
This violates property (e) of Brownian web and its dual listed earlier.
Hence $\P(\Omega^{k}_{0}) = 0$ for all $k \geq 0$ and this completes the proof that
$\eta_n \to \eta$ almost surely as $n \to \infty$.

Now it suffices to show that the sequence $\{\eta_n : n \in \N\}$ is uniformly integrable. As $\eta_n \to \eta$ almost surely, using Proposition 4.1 of \cite{FINR04}
 for any $k \geq 1$ we have,
\begin{align*}
\lim_{n \to \infty}\P(\eta_n \geq k) = \P(\eta \geq k) \leq (\P(\eta \geq 1))^k < 1.
\end{align*}
We have used the fact that $\P(\eta \geq 1)$ is same as the probability that two independent Brownian motions starting at unit distance do not meet by time $1$, which is strictly smaller 
than $1$. This completes the proof.
\qed

\begin{remark}
\label{rem:UIRobustArgument}
It is important to observe that the argument for uniform integrability presented 
here depends only on the fact that $\eta_n \to \eta$ almost surely whereas the argument of Lemma 3.3 of \cite{RSS15} uses some model specific assumptions.
We comment here that similar arguments show 
that $\eta_n(t_0,t_1) \to \eta(t_0,t_1)$ almost surely for any $t_0 < t_1$.
\end{remark}

The following lemma proves Theorem \ref{thm:ScalingLimCRTJt}.
\begin{lemma}
\label{lem:MetricExpConv}
As $n \to \infty$, we have  $\E[\kappa_n(f)] \to \E[\kappa(f)] $.
\end{lemma}
We first prove Theorem \ref{thm:ScalingLimCRTJt} assuming Lemma \ref{lem:MetricExpConv} holds.

\noindent \textbf{Proof of Theorem \ref{thm:ScalingLimCRTJt}}:
Recall that for $(x,t) \in \Z^2_{\text{even}}$ and for $\pi^{(x,t)} \in {\cal X}$, i.e.,
the path in ${\cal X}$ starting from $(x,t)$, the $n$ -th order diffusively scaled path 
is denoted by $\pi^{(x,t)}_n$. For $(x,t) \in \Z^2_{\text{even}}$ and for $n \geq 1$, let 
\begin{align*}
S_n(x,t) := \{& \pi^{(y,s)}_n \in {\cal X}_n : (y,s) \in C(x,t) \},\text{ where }C(x,t)\text{ is defined as in (\ref{eq:ClusterLevel}) and }\\
\widehat{S}_n(x,t) := \{ & \widehat{\pi}_n : \text{ where }\widehat{\pi}_n \text{ is the }n\text{-th order
scaled version of }\widehat{\pi}\in  \widehat{\cal X}\text{ such that } \\
&\widehat{\pi}^{(x-1,t)}(s)\leq \widehat{\pi}(s) \leq \widehat{\pi}^{(x+1,t)}(s)\text{ for all } s\leq \sigma_{\widehat{\pi}}  \}.
\end{align*}
Both the path spaces, $S_n(x,t)$ and $\widehat{S}_n(x,t)$, are tree-like in the sense of 
Definition \ref{def:TreeLikeMetricFrmCoalPaths} and we consider the complete metric spaces 
${\cal M}(S_n(x,t))$ and ${\cal M}(\widehat{S}_n(x,t))$. 
It follows that for all $(x,t) \in \Z^2_{\text{even}}$ and for all $n \geq 1$, we have
that both the metric spaces $T_n(x,t)$ and $\widehat{T}_n(x,t)$ are naturally embedded into 
 the metric spaces ${\cal M}(S_n(x,t))$ and ${\cal M}(\widehat{S}_n(x,t))$ respectively.
 This embedding ensures that for all $n \geq 1$ we have 
$$
d_{\text{GH}}\bigl({\cal M}(S_n(x,t)),T_n(x,t)\bigr)\vee d_{\text{GH}}\bigl({\cal M}(\widehat{S}_n(x,t)),\widehat{T}_n(x,t)\bigr)
\leq 1/n.
$$
Hence to prove Theorem \ref{thm:ScalingLimCRTJt} it suffices to show that 
$$
({\cal M}(S_n(0,0)),{\cal M}(\widehat{S}_n(0,0)))|\{L(0,0) \geq n\}\Rightarrow (\TT, \widehat{\TT})\text{ as }n \to \infty.
$$ 
Using translation invariance of our model, we have
\begin{align*}
\E[\kappa_n(f)] = & \E\Bigl[\sum_{k = 0}^{\lfloor \sqrt{n} \rfloor} \mathbf{1}_{\{L(k,0) \geq n\}}f\bigl(\phi_k\bigl({\cal M}(S_n(k,0))\bigr), \phi^\prime_k \bigl({\cal M}(\widehat{S}_n(k,0))\bigr)\bigr)\Bigr]\\
& \qquad \qquad \text{ where }\phi_k,\phi^\prime_k \text{ are isometric embeddings of the respective metric spaces into }\mathbb{M}\\
= & (\lfloor \sqrt{n}\rfloor +1) \E\Bigl[\mathbf{1}_{\{L(0,0) \geq n\}}f\bigl(\phi_0\bigl({\cal M}(S_n(0,0))\bigr), \phi^\prime_0 \bigl({\cal M}(\widehat{S}_n(0,0))\bigr)\bigr)\Bigr]\\
= & (\lfloor \sqrt{n}\rfloor +1)\P(L(0,0)\geq n)\E\Bigl[f\bigl(\phi_0\bigl({\cal M}(S_n(0,0))\bigr), \phi^\prime_0 \bigl({\cal M}(\widehat{S}_n(0,0))\bigr)\bigr)\Bigl|L(0,0)\geq n\Bigr]\\
= & \E[\sum_{k = 0}^{\lfloor \sqrt{n} \rfloor} \mathbf{1}_{\{L(k,0) \geq n\}}]
\E\Bigl[f\bigl(\phi_0\bigl({\cal M}(S_n(0,0))\bigr), \phi^\prime_0 \bigl({\cal M}(\widehat{S}_n(0,0))\bigr)\bigr)\Bigl|L(0,0)\geq n\Bigr]\\
= & \E[\eta_n]\E\Bigl[f\bigl(\phi_0\bigl({\cal M}(S_n(0,0))\bigr), \phi^\prime_0 \bigl({\cal M}(\widehat{S}_n(0,0))\bigr)\bigr)\Bigl|L(0,0)\geq n\Bigr].
\end{align*}
Using Lemma \ref{lem:OccupdPtSetExpConv} and Lemma \ref{lem:MetricExpConv}
we obtain 
\begin{align*}
\lim_{n\to \infty} \E\Bigl[f\bigl(\phi_0\bigl({\cal M}(S_n(0,0))\bigr), \phi^\prime_0 \bigl({\cal M}(\widehat{S}_n(0,0))\bigr)\bigr)\Bigl|L(0,0)\geq n\Bigr] & = 
\lim_{n\to \infty} \frac{ \E[\kappa_n(f)] }{ \E[\eta_n] } \\
\\ & = \frac{ \E[\kappa(f)] } { \E[\eta] } = \E ( f(\phi(\TT),\phi^\prime(\widehat{\TT}) )).
\end{align*}
Since $f$ is chosen arbitrarily, this completes the proof.\qed

%

We now prove Lemma \ref{lem:MetricExpConv}. 
Since $f$ is bounded, uniform integrability of $\{\kappa_n(f)\}$ follows from uniform integrability of $\{\eta_n: n \in \N\}$.
Hence to prove Lemma \ref{lem:MetricExpConv},
it suffices to show that $\kappa_n(f) \to \kappa(f)$ almost surely as $n \to \infty$.   
We introduce some notation now.
On the event $\{\eta = 1\}$, Lemma \ref{lem:OccupdPtSetExpConv}  ensures that there exists $n_0 = n_0(\omega)$ such that $\eta_n = 1$ for all $n\geq n_0$.
For $n \geq n_0$ consider $x=x(\omega), x_n = x_n(\omega) \in \R$
such that $\xi = \{x\}$ and $\xi_n = \{ x_n \}$.
For ease of notation on the event $\{\eta = 1\}$ and for $n\geq n_0$ we take
\begin{align*}
(S, \widehat{S}) & := \bigl({\cal S}(\widehat{\pi}^{(x,0)}_r,\widehat{\pi}^{(x,0)}_l)_{(\WW, \widehat{\WW})},
\widehat{\cal S}(\widehat{\pi}^{(x,0)}_r,\widehat{\pi}^{(x,0)}_l)_{(\WW, \widehat{\WW})}\bigr),\\
(\overline{S}, \overline{\widehat{S}}) & := \bigl(\overline{{\cal S}}(\widehat{\pi}^{(x,0)}_r,\widehat{\pi}^{(x,0)}_l)_{(\WW, \widehat{\WW})},\overline{\widehat{\cal S}}(\widehat{\pi}^{(x,0)}_r,\widehat{\pi}^{(x,0)}_l)_{(\WW, \widehat{\WW})}\bigr),\\ 
(\overline{S}_n, \overline{\widehat{S}}_n) = (S_n, \widehat{S}_n) & := 
\bigl({\cal S}(\widehat{\pi}^{(x_n,0)}_r,\widehat{\pi}^{(x_n,0)}_l)_{({\cal X}_n, 
\widehat{{\cal X}}_n)}, \widehat{{\cal S}}(\widehat{\pi}^{(x_n,0)}_r,\widehat{\pi}^{(x_n,0)}_l)_{({\cal X}_n, \widehat{{\cal X}}_n)} \bigr),\\
({\cal M}, \widehat{{\cal M}}) & := ({\cal M}^{(x,0)}_{(\WW, \widehat{\WW})}, 
\widehat{{\cal M}}^{(x,0)}_{(\WW,\widehat{\WW})})\text{ and }\\
({\cal M}_n, \widehat{{\cal M}}_n) & := ({\cal M}^{(x_n,0)}_{({\cal X}_n, 
\widehat{{\cal X}}_n)}, \widehat{{\cal M}}^{(x_n,0)}_{({\cal X}_n, \widehat{{\cal X}}_n)}).\end{align*} 
In what follows, we assume that we are working on 
a probability space such that $(\bar{{\mathcal  X}}_n, \bar{\widehat{{\mathcal  X}}}_n)$ converges  
to $({\mathcal  W}, \widehat{{\mathcal  W}})$ almost surely in $({\mathcal  H}\times \widehat{\mathcal  H},d_{{\mathcal  H}\times \widehat{{\mathcal  H}}})$.
\begin{lemma}
\label{lem:ConvSnS}
On the event $\{\eta = 1\}$, we have
$d_{{\mathcal H}\times \widehat{{\mathcal H}}}((S_n, \widehat{S}_n), (\overline{S}, \overline{\widehat{S}})) \to 0$
almost surely  as $n \to \infty$.
\end{lemma}
\noindent \textbf{ Proof :} 
We prove that $d_{{\mathcal  H}}(S_n, \overline{S}) \to 0$ as $n \to \infty$
almost surely. The argument for $d_{\widehat{\mathcal  H}}(\widehat{S}_n, \overline{\widehat{S}})$ is exactly the same and hence omitted. 
Fix $ \epsilon > 0 $. 
Choose $\delta =  \delta(\omega) \in (0, \epsilon)\cap \Q $
such that 
$$\sup\{\widehat{\pi}(s_1) - \widehat{\pi}(s_2): 
s_1,s_2 \in [- 2\delta, 0], \widehat{\pi}\in \widehat{\WW}^{0+}, \widehat{\pi}(0)
\in [0,1] \} < \epsilon/8.$$
Assume that $ \xi(-\delta, 0) = \{ x, x_1,  \dotsc, x_{k} \}$ where $\xi = \{x\}$.
We choose $ \gamma_{\delta} = \gamma_{\delta}(\omega) \in (0,\epsilon)$ 
such that (a) the points in $\xi(-\delta, 0)$ are at least $2\gamma_\delta$ distance apart, i.e., $  | x - x_i |\wedge |x_i - x_j| > 2 \gamma_{\delta} $ for all $1\leq i < j \leq k$  and  (b) $(x- \gamma_{\delta}, x +\gamma_{\delta})$ as well as $(x_i- \gamma_{\delta},
x_i +\gamma_{\delta})$  are both included in the interval $(0,1)$ for all $1\leq i \leq k$.
We mentioned earlier that there exist 
dual paths  $\widehat{\pi}^{(x, 0)}_r,\widehat{\pi}^{(x,0)}_l$ in $\widehat{\WW}$
both starting at $ (x, 0)$ such that 
 $\widehat{\pi}^{(x, 0)}_r(-1) >  
\widehat{\pi}^{(x,0)}_l(-1) $. Similarly for each $x_i \in \xi(-\delta, 0)$ there exist dual paths  $\widehat{\pi}^{(x_i, 0)}_r, 
\widehat{\pi}^{(x_i,0)}_l$  in $\widehat{\WW}$ both starting at $ (x_i, 0)$ such that 
$\widehat{\pi}^{(x_i, 0)}_r(-\delta) >  \widehat{\pi}^{(x_i, 0)}_l(-\delta) $ 
 for all $ 1 \leq i \leq k$. Because of continuity, there exists $\nu_{\delta} = \nu_{\delta}(\omega) > 0
 $ such that $\widehat{\pi}^{(x,0)}_r( -1 - \nu_{\delta}) >  \widehat{\pi}^{(x, 0)}_l(- 1- \nu_{\delta}) $ and $\widehat{\pi}^{(x_i,0)}_r( -\delta - \nu_{\delta}) >  \widehat{\pi}^{(x_i, 0)}_l(- \delta- \nu_{\delta}) $ for $1 \leq i \leq k$. 
 From properties (f) and (g) of $(\WW, \widehat{\WW})$ (see Subsection \ref{subsec:Bweb}) it follows that there exist $\widehat{\pi}_{l},\widehat{\pi}_{r}, \widehat{\pi}^{i}_{l},\widehat{\pi}^{i}_{r} \in \widehat{\cal W}$ such that all these dual paths have starting times strictly 
larger than $0$ and 
\begin{align}
\label{eq:Calculation3}
\begin{split}
& \text{(i) }\max\{|(\widehat{\pi}_{l} -\widehat{\pi}^{(x, 0)}_l) (0)|, |(\widehat{\pi}_{r} -\widehat{\pi}^{(x, 0)}_r) (0)|, |(\widehat{\pi}^{i}_{l} - \widehat{\pi}^{(x_i, 0)}_l) (0)|,
|(\widehat{\pi}^{i}_{r} - \widehat{\pi}^{(x_i, 0)}_r )(0)|:1\leq i \leq k \} < \gamma_\delta /4,\\
& \text{(ii) }\max\{|(\widehat{\pi}_{l} -\widehat{\pi}^{(x, 0)}_l) (-\delta)|, |(\widehat{\pi}_{r} -\widehat{\pi}^{(x, 0)}_r) (-\delta)|, |(\widehat{\pi}^{i}_{l} - \widehat{\pi}^{(x_i, 0)}_l) (-\delta)|, |(\widehat{\pi}^{i}_{r} - \widehat{\pi}^{(x_i, 0)}_r )(-\delta)|:1\leq i \leq k \} = 0.
\end{split}
\end{align}
 This gives that 
\begin{align}
\label{eq:Calculation4}
\max\{(\widehat{\pi}_{r} - \widehat{\pi}_{l})(s): s \in [-2\delta, 0]\}
\vee \max\{(\widehat{\pi}^{i}_{r} - \widehat{\pi}^{i}_{l})(s): s \in [-2\delta, 0],
 1\leq i\leq k\} < \epsilon/2 .
\end{align}
From Remark \ref{rem:UIRobustArgument} we have $\P\bigl(\bigcap_{t \in [-1,0)\cap \Q}\{\lim_{n \to \infty}\eta_n(t,0) = \eta(t,0)\}\bigr) = 1$.  
Set $\zeta_{\delta} > 0$ such that
\begin{align*} 
\zeta_{\delta}  < \min\{ (\widehat{\pi}^{(x,0)}_r -\widehat{\pi}^{(x , 0)}_l)
(-1 -\nu_{\delta}), \sigma_{\widehat{\pi}_{l}}, \sigma_{\widehat{\pi}_{r}}\}\wedge 
\min\{ (\widehat{\pi}^{(x_i, 0)}_r -
\widehat{\pi}^{(x_i, 0)}_l)(-\delta - \nu_{\delta}),
\sigma_{\widehat{\pi}^{i}_{l}}, \sigma_{\widehat{\pi}^{i}_{r}} : 1 \leq i \leq k \}.
\end{align*}
Let $n_1 = n_1(\omega) >  \max\{4/\nu_{\delta},1/\delta, n_0\}$ be such that, 
for all $n \geq n_1$
\begin{itemize}
 \item [(i)] $\eta_n(- \delta,0) = \eta(- \delta,0)$ and $\eta_n = \eta$;
 
 \item[(ii)] $d_{{\cal H}\times \widehat{\cal H}}(({\cal W},\widehat{{\cal W}}),(\bar{{\cal X}}_n,
\bar{\widehat{{\cal X}}}_n)) < g(\min \{ \gamma_{\delta}, \zeta_{\delta}, \epsilon \}/4,m)$ where $m$ is defined as in (\ref{eq:dualCompact}) and the value $g(\min \{ \gamma_{\delta}, \zeta_{\delta}, \epsilon \}/4,m)$ is taken as in Remark \ref{rem:MetricConv}.
\end{itemize}

The choice of $n_1 $ and $ \zeta_{\delta}$ ensure that 
there exist $\widehat{\pi}^n_{r}, \widehat{\pi}^n_{l} \in \widehat{\cal X}^{0+}_n$, and 
for all $1 \leq i \leq k$ there exist 
$\widehat{\pi}^n_{i,r}, \widehat{\pi}^n_{i,l} \in \widehat{\cal X}^{0+}_n$ such that
\begin{align}
\label{eqn:DualApproxDist} 
\sup\{|(\widehat{\pi}^n_{i,r} - \widehat{\pi}^{i}_{r})(s)|, 
|(\widehat{\pi}^n_{i,l} - \widehat{\pi}^{i}_{l})(s)|,
|(\widehat{\pi}^n_{r} - \widehat{\pi}_{r})(s)|, 
|(\widehat{\pi}^n_{l} - \widehat{\pi}_{l})(s)|: s \in [-1,0] \}< \min \{ \gamma_{\delta}, \zeta_{\delta}, \epsilon\}/4,  
\end{align}
as well as
$\widehat{\pi}^n_{r}(-1 ) > \widehat{\pi}^n_{l}(-1)$ and 
$\widehat{\pi}^n_{i,r}( -\delta -\nu_{\delta}) > \widehat{\pi}^n_{i,l}( -\delta -\nu_{\delta})$ for all $ 1 \leq i \leq k$. 
Therefore,  there exists a forward path $ \vartheta^n \in {\cal X}^{(-1)-}_n $ such that $ \widehat{\pi}^n_{l}(-1) < \vartheta^n (-1) < \widehat{\pi}^n_{r}(-1)  $. Using (\ref{eq:Calculation3}) we further have that $x - \gamma_{\delta} < \widehat{\pi}^n_{l}(0) < \vartheta^n (0) < 
\widehat{\pi}^n_{r}(0) < x + \gamma_{\delta} $. 
Similarly for each $1\leq i \leq k$ there exists a forward path $ \vartheta^n_i \in {\cal X}^{(-\delta)-}_n $ such that (a) $ \widehat{\pi}^n_{i,l}(-\delta) < \vartheta^n_i (-\delta) < \widehat{\pi}^n_{i,r}(-\delta)  $
and (b) $x_i - \gamma_{\delta} < \widehat{\pi}^n_{i,l}(0) < \vartheta^n_i (0) < 
\widehat{\pi}^n_{i, r}(0) < x_i + \gamma_{\delta} $.
Since for all $n \geq n_1$ we have $\# \xi_n(-\delta, 0) = k + 1$,
hence $ \xi_n(-\delta, 0) = \{ \vartheta_i^n (0) : 1 \leq i \leq k  \} \cup\{\vartheta^n (0)\} 
$. Thus, 
\begin{align}
\label{x_n_unique}
& \text{ there exist a unique $x^n = \vartheta^n (0) \in 
\xi_n(-\delta, 0)$ with $|x^n - x| < \gamma_{\delta}$ and } \nonumber\\
& \text{ unique }x_i^n = \vartheta_i^n (0) \in 
\xi_n(-\delta, 0) \text{ with $|x_i^n - x^i| < \gamma_{\delta}$ for each $i$. } 
\end{align}

We observe that $(y^n,0) = (\sqrt{n}\vartheta^n (0),0)$ is in $\Z^2_{\text{even}}$
and the vertices $(y^n + 1, 0), (y^n - 1, 0)$ are the right and left dual neighbours of 
$(y^n,0)$ respectively. Similarly $(y^n_i,0) = (\sqrt{n}\vartheta^n_i (0),0) \in \Z^2_{\text{even}}$ for $1\leq i \leq k$.

For any $\pi^n \in S_n$ we have $\pi^n(0) = x^n$ and hence for all $n \geq n_1$, using (\ref{eq:Calculation4}) and (\ref{eqn:DualApproxDist}) obtain
\begin{align*}
\sup_{ s_1, s_2 \in [(-2\delta)\vee\sigma_{{\pi}^n}, 0]}\{|{\pi}^n(s_1)-{\pi}^n(s_2)|\} & \leq
\sup_{ s_1, s_2 \in [-2\delta, 0]}\{|\widehat{\pi}_n^{(y^n + 1,0)}(s_1) - \widehat{\pi}_n^{(y^n - 1,0)}(s_2)|\} \\
& \leq \sup_{ s_1, s_2 \in [-2\delta, 0]}\{|\widehat{\pi}^n_r(s_1) - \widehat{\pi}^n_l(s_2)|\}\\
& \leq\sup_{ s_1, s_2 \in [-2\delta, 0]}\{|\widehat{\pi}_r(s_1) - \widehat{\pi}_l(s_2)|\} + \epsilon/2 \leq \epsilon.
\end{align*}
For any $\pi_1, \pi_2 \in \overline{S}$ and $\pi^n_1, \pi^n_2 \in S_n$ we have 
$\pi_1(s) = \pi_2(s)$ and $\pi^n_1(s) = \pi^n_2(s)$ for all $s \geq 0$ with $|\pi^n_1(0) - \pi_1(0)|< \epsilon$. This observation together with the fact that $\sup_{ s_1, s_2 \in [(-2\delta)\vee\sigma_{{\pi}_1}, 0]}\{|{\pi}_1(s_1)-{\pi}_1(s_2)|\}\leq \sup_{ s_1, s_2 \in [-2\delta, 0]}\{|\widehat{\pi}_r(s_1)-\widehat{\pi}_l(s_2)|\}< \epsilon$ ensures that
to prove Lemma \ref{lem:ConvSnS}, it suffices to show that
$$
d_{{\cal H}}(S^{(-2\delta)-}_n, \overline{S}^{(-2\delta)-}) \to 0\text{ as }n\to \infty.
$$

Fix $\pi \in \overline{S}$ with $\sigma_\pi \leq - 2\delta$. 
From the choice of $n_1$ it follows that for all $n \geq n_1$ there exists  an 
 approximating path $\pi^n \in {\cal X}^{(-\delta)-}_n$ with $|\pi(0) - \pi^n(0)| 
 < \gamma_{\delta}$. As $\sigma_{\pi^n}\leq -\delta$, we have $\pi^n(0) \in \xi(-\delta, 0)$. 
 Since $\pi \in \overline{S}$, it 
follows that $\pi(0) = x$. Hence from the uniqueness as mentioned in (\ref{x_n_unique}), it follows that for all $n \geq n_1$ the approximating path $\pi^n$ must have 
$\pi^n(0) = x^n$ implying that $\pi^n \in S_n$. 
Similar argument shows that for all $n \geq n_1$ and for $\pi^n \in S_n$ with $\sigma_{\pi^n} \leq - \delta$, there exists approximating path in $\overline{S}$.  
This completes the proof.
\qed

Before stating the next lemma we give an alternate definition of $d_{GH}$ which will be used in our proof. For two compact metric spaces $(X_1 , d_1 )$ and $(X_2 , d_2 )$, a correspondence between
$X_1$ and $X_2$ is a subset ${\cal R}$ of $X_1 \times X_2$ such that for every $x_1 \in X_1$ there exists at least one $x_2 \in X_2$ such that $(x_1 , x_2 ) \in {\cal R}$ and conversely for every $y_2 \in X_2$
there exists at least one $y_1 \in X_1$ such that $(y_1 , y_2 ) \in {\cal R}$. The distortion of a correspondence ${\cal R}$ is defined by
$$
\text{dis}({\cal R}) := \sup\{|d_1(x_1 , y_1) -  d_2(x_2 , y_2 )| : (x_1 , x_2 ), (y_1 , y_2 ) \in {\cal R}\}.
$$
Then, if $X_1$ and $X_2$ are two real trees, we have
\begin{align}
\label{eq:GHDistCorr}
d_{\text{GH}} (X_1 , X_2 ) := 1/2 \inf_{
{\cal R}\in C(X_1 ,X_2 )} \text{dis}({\cal R})
\end{align}
where $C(X_1 , X_2 )$ denotes the set of all correspondences between $X_1$ and $X_2$.
Using Lemma \ref{lem:ConvSnS} we prove the
following lemma which is the main lemma of this section. 
\begin{lemma}
\label{lem:GHConvMnM}
On the event $\{\eta = 1\}$ we have 
$d_{\text{GH}}(\phi_n({\cal M}_n),\phi({\cal M}))\vee 
d_{\text{GH}}(\phi^\prime_n(\widehat{\cal M}_n),\phi^\prime(\widehat{\cal M}))\to 0$ 
almost surely  as $n \to \infty$, where $\phi_n,\phi,\phi^\prime_n, \phi^\prime$ denote the corresponding embeddings of the respective metric spaces into $\mathbb{M}$.
\end{lemma}
\noindent \textbf{ Proof :} Before proving this lemma we first sketch the key idea.
The metric $d_{{\mathcal H}\times \widehat{{\mathcal H}}}$ deals with the path space topology, convergence in this metric does not ensure that the coalescing times also converge (even if they are finite). Using joint convergence of the collection of forward paths together with the collection of dual paths, we 
 show that the coalescing times of the scaled paths converge jointly to the coalescing time of the limiting Brownian paths and this gives convergence with respect to the ancestor metric. We first show that $d_{\text{GH}}(\phi_n({\cal M}_n), \phi({\cal M})) \to 0$ as $n \to \infty$ almost surely and the argument for
$d_{\text{GH}}(\phi^\prime_n(\widehat{{\cal M}}_n), \phi^\prime(\widehat{\cal M})) \to 0$ as $n \to \infty$ is exactly the same.

Using Definition (\ref{eq:GHDistCorr}), it suffices to show that 
there exists a correspondence ${\cal R}_n $ between 
${\cal M}_n$ and ${\cal M}$ such that
\begin{align}
\label{eq:DistCorrLim}
\lim_{n \to \infty} \text{dis}({\cal R}_n) = 0.
\end{align}

Choose any sequence $\{\epsilon_k : k \in \N\}$ such that
$\epsilon_k \downarrow 0$ as $k \to \infty$. For $k \geq 1$, let $n_k = n_k(\omega) \in \N$ be such that $d_{\cal H}(S_{n_k} , \overline{S}) < g(\epsilon_k/8,m)$ almost surely where 
$m$ is defined as in (\ref{eq:dualCompact}) and $g(\epsilon,m)$ as in Remark \ref{rem:MetricConv}. 
We define correspondences  ${\cal R}_n$ for large $n$ only.
For $n_k \leq n < n_{k+1} $, we define our correspondence ${\cal R}_n$ between ${\cal M}_{n}$
 and ${\cal M}$ as
\begin{align}
\label{eq:k-Corr}
{\cal R}_n := & \Bigl(\bigcup_{(y,s) \in {\cal M}}\{((y,s),(y_{n},s_{n})): 
 (y_{n},s_{n}) \in {\cal M}_{n}\text{ such that }
d_{\Pi}(\pi^{(y^\prime, s^\prime)},\pi^{(y_{n},s_{n})}) < g(\epsilon_k/8,m)\nonumber\\
& \qquad \qquad \text{ for some }
\pi^{(y^\prime, s^\prime)} \in S \text{ with }d_A((y, s),(y^\prime, s^\prime))<\epsilon_k/8\}\Bigr)\bigcup\nonumber\\
& \Bigl(\bigcup_{(y_{n},s_{n}) \in {\cal M}_{n}}\{((y,s),(y_{n},s_{n})): 
 (y,s) \in {\cal M}\text{ with }\pi^{(y, s)} \in S \text{ such that }\nonumber\\
& \qquad \qquad d_{\Pi}(\pi^{(y, s)},\pi^{(y_{n},s_{n})})< g(\epsilon_k/8,m)\} \Bigr) . 
\end{align}
For all $n_k \leq n < n_{k+1}$, (\ref{eq:k-Corr}) gives a valid correspondence.

Aiming for a contradiction suppose (\ref{eq:DistCorrLim}) does not hold. 
As $d_{\Pi}(\pi_1, \pi_2) < g(\epsilon_{k}/8,m)$ implies that 
 $|\sigma_{\pi_1}- \sigma_{\pi_2}| < \epsilon_{k}/8 $, 
if (\ref{eq:DistCorrLim}) does not hold then there exists
 $\beta = \beta(\omega) >0$ such that  
\begin{align}
\label{eq:InfBeta}
\limsup_{k \to \infty}
\sup\{ & |\gamma(1,2)_k - \gamma(1,2)_{n}| : \pi^k_1, \pi^k_2 \in S\text{ and } \pi^{n}_1, \pi^{n}_2 \in S_{n} \text{ for some }n_k \leq n < n_{k+1} \text{ and }\nonumber\\
 & \qquad d_{\Pi}(\pi_1, \pi^{n}_1)\vee d_{\Pi}(\pi_2, \pi^{n}_2) < g(\epsilon_k/8, m)\} > \beta, 
\end{align}
where for ease of notation we take $\gamma(1,2)_k = \gamma(\pi^k_1,\pi^k_2)$ and $\gamma(1,2)_{n} = \gamma(\pi^{n}_1,\pi^{n}_2)$.
Choose $k_0$ such that $\epsilon_k \vee (1/n_{k})< \beta/8$ for all $k \geq k_0$.
From (\ref{eq:InfBeta}), it follows that there are infinitely many $k \geq k_0$ such that there exist $\pi^{k}_1, \pi^{k}_2 \in S$ and $\pi^{n}_1, \pi^{n}_2 
 \in S_{n}$ for some $n_k \leq n < n_{k+1}$ with $d_{\Pi}(\pi^{k}_1, \pi^{n}_1)\vee d_{\Pi}(\pi^{k}_2, \pi^{n}_2) < g(\epsilon_{k}/8,m)$ and $|\gamma(1,2)_{k} - \gamma(1,2)_{n}| > \beta/2$.
We consider the two cases :
\begin{itemize}
 \item[(a)] $\gamma(1,2)_{n} > \gamma(1,2)_{k} + \beta$  and 
 \item[(b)] $\gamma(1,2)_{n} \leq \gamma(1,2)_{k} - \beta$.
\end{itemize}

For (a) from the choice of $k_0$, it follows that for infinitely many $k \geq k_0$
and for some $n_k \leq n < n_{k+1}$
we have $\sigma_{\pi^{n}_1}\vee \sigma_{\pi^{n}_2} \leq \gamma(1,2)_k + \beta/4$
with $\pi^{n}_1(\gamma(1,2)_k + \beta) \neq \pi^{n}_2(\gamma(1,2)_k + \beta)$.
From the construction of the dual graph it follows that for infinitely many $k\geq k_0$
and for some $n_k \leq n < n_{k+1}$, 
there exists $\widehat{\pi}^{n} \in \widehat{S}_{n}$ with $\sigma_{\widehat{\pi}^{n}} \geq \gamma(1,2)_k + 3\beta/4$ and for $s \in [\gamma(1,2)_k, \gamma(1,2)_k + 3\beta/4]$ we have $\widehat{\pi}^{n}(s) \in (\pi^{n}_1(s)\wedge \pi^{n}_2(s), \pi^{n}_1(s)\vee \pi^{n}_2(s))$.
 Hence we have  
$|\widehat{\pi}^{n}(s) - \pi^{k}_1(s)| < \epsilon_k/2$ for all 
$s \in [\gamma(1,2)_k, \gamma(1,2)_k + 3\beta/4]$.
Since $d_{{\mathcal H}\times \widehat{{\mathcal H}}}((S_{n},\widehat{S}_{n}),(\overline{S}, \overline{\widehat{S}})) < g(\epsilon_{k}/8,m) $ for all $n \geq n_k$, 
for infinitely many $k \geq k_0$ there exists $\widehat{\pi}^{k} \in \widehat{S}$ with $\sigma_{\widehat{\pi}^{k}} \geq \gamma(1,2)_k + \beta/4$
and $d_{\widehat{\Pi} }(\widehat{\pi}^{k}, \widehat{\pi}^{n_k}) < \epsilon_{k}/8$.
This gives us that $|\widehat{\pi}^{k}(s) - \pi^{k}_1(s)| < \epsilon_{k}$ for all 
$s \in [\gamma(1,2)_k, \gamma(1,2)_k + \beta/4]$.

For (b) we need a slightly different argument.
From the properties of $(\WW, \widehat{\WW})$ discussed earlier, it follows that there exists
$\widehat{\pi}^{k} \in \widehat{S}$ starting from the point $(\pi^{k}_1(\gamma(1,2)_k), \gamma(1,2)_k)$ and lying  between the forward paths $\pi^{k}_1$ and $\pi^{k}_2$ in the time interval $[\sigma_{\pi^{k}_1} \vee \sigma_{\pi^{k}_2}, \gamma(1,2)_k]$. 
For $s \in [\gamma(1,2)_k - \beta/4,\gamma(1,2)_k]$ we have
 $\pi^{n}_1(s) = \pi^{n}_2(s)$. Since
$d_{\Pi}(\pi^{k}_i, \pi^{n}_i) < g(\epsilon_{k}/8,m)$ for $i=1,2$,
it follows that $0<|\pi^{k}_1(s) - \pi^{k}_2(s)| < \epsilon_k/4$ for all 
$s \in [ \gamma(1,2)_k - \beta/4, \gamma(1,2)_k]$.
Since the dual path $\widehat{\pi}^{k}$ starting from the point $(\pi^{k}_1(\gamma(1,2)_k), \gamma(1,2)_k)$ lies between $\pi^{k}_1$ and $\pi^{k}_2$, we have that 
$|\widehat{\pi}^{k}(s) - \pi^{k}_1(s)| < \epsilon_{k}/4$ for all 
$s \in [ \gamma(1,2)_k - \beta/4, \gamma(1,2)_k]$.

Since $\epsilon_k \downarrow 0$, because of compactness of $(\WW, \widehat{\WW})$, 
this contradicts the fact that in the double Brownian web almost surely no two 
forward path and dual path spend positive Lebesgue measure time together.
This completes the proof.
\qed

Finally we prove Lemma \ref{lem:MetricExpConv} to complete the proof of Theorem
\ref{thm:ScalingLimCRTJt}.

\noindent \textbf{ Proof of Lemma \ref{lem:MetricExpConv} :}
Since $f$ is bounded, the uniform integrability of the sequence $\{\kappa_n(f) ; n \in \N\}$
follows from that of the sequence $\{\eta_n : n \in \N\}$. 
Since $\eta_n \to \eta$ almost surely, $\kappa_n(f) \to \kappa(f)$ holds trivially on the event $\{\eta = 0\}$.
On the event $\{\eta = 1\}$, using continuity of $f$ this follows from 
Lemma \ref{lem:GHConvMnM}). For general $k \geq 1$ on the event $\{\eta = k\}$, the proof is similar.
%
\qed

\section{Concluding remark}
\label{sec:Conclusion}

To the best of our knowledge 
both these continuum random trees, $\TT$ and $\widehat{\TT}$, have not been studied
in the literature so far. On the other hand, the system of coalescing Brownian motions starting from every space-time points on $\mathbb{R}^2$ has been 
extensively studied but with respect to a different topology (see 
\cite{FINR04}, \cite{SSS16} for details).
It should be mentioned here that coalescing Brownian motions starting from all the points in $\R$
at a given time has been studied as a genealogical tree in \cite{GSW16} with respect
to a different topology where the spatial locations of the points are also taken care of. 

Further understanding of these continuum random trees $\TT$ and $\widehat{\TT}$
might be useful and we present some questions towards that direction.
In the above construction we have used the fact that the collection of 
 forward paths ${\cal S}(\widehat{B}^+,\widehat{B}^-)$ in the region $\Delta(\widehat{B}^+,\widehat{B}^-)$ is almost surely uniquely determined by the 
 collection of backward or dual paths $\widehat{{\cal S}}(\widehat{B}^+,\widehat{B}^-)$.
It should be possible to construct $\TT$ 
 directly as completion of the metric space obtained from the system of 
 coalescing forward Brownian paths starting from all the points of 
 $\Delta(\widehat{B}^+,\widehat{B}^-)\cap \Q^2$,
 which follow Skorohod reflection at the boundary 
 of $\Delta(\widehat{B}^+,\widehat{B}^-)$. In order to do that construction, one 
 has to show that starting from finitely many space-time points such a collection satisfies
 Kolmogorov's consistency conditions, as it does \textit{not} directly follow from \cite{STW00}.

 From the construction of $(\TT, \widehat{\TT})$, it is reasonable to expect that $\widehat{\TT}$ almost surely determines $\TT$ and vice-versa. One possible to way to prove this is to show that 
 the contour function of $\TT$ is determined by $\widehat{\TT}$ and vice-versa and in that case $\widehat{\TT}$ can be regarded as the \textit{dual} tree of $\TT$.
 It is natural to ask how are the contour functions of these two continuum random trees distributed?  Presently we {\cg \textit{do not}} have any understanding about these processes.   
Finally in the context of drainage network scaling relations, it might be of interest to
 see whether these two continuum trees obey Horton's law or not.

{\bf Acknowledgements:} 
The author thanks Siva Athreya, Manjunath Krishnapur, Rahul Roy, Rongfeng Sun and Sreekar Vadlamani for discussion and comments. Part of the work was done when the author was a visiting scientist at Indian Statistical Institute, Bangalore center.

\end{document}